\documentclass[12pt]{article}

\pdfoutput=1

\usepackage{mathtools}
\usepackage{booktabs}
\usepackage[english]{babel} 
\usepackage[protrusion=true,expansion=true]{microtype} 
\usepackage{amsmath,amsfonts,amsthm}
\usepackage{ amssymb }
\usepackage{graphicx}
\usepackage{array}
\usepackage{bbm}
\usepackage{multirow}
\usepackage{hhline}
\usepackage[titletoc]{appendix}

\usepackage{booktabs} 
\usepackage{natbib}
\usepackage{setspace}
\usepackage{graphicx}
\usepackage{subcaption}

\usepackage{microtype} 
\usepackage{tabularx}
\usepackage[font = small,labelfont=bf,textfont=it]{caption} 
\usepackage{footnote}
\usepackage{algorithm}
\usepackage[noend]{algpseudocode}

\makeatletter
\newcommand{\distas}[1]{\mathbin{\overset{#1}{\kern\z@\sim}}}%
\newsavebox{\mybox}\newsavebox{\mysim}
\newtheorem{theorem}{Theorem}[section]

\newtheorem{lemma}[theorem]{Lemma}
\theoremstyle{definition}

\newtheorem{assumption}{Assumption}[section]

\newtheorem{proposition}{Proposition}[section]

\newtheorem{remark}[theorem]{Remark}
\newcommand{\distras}[1]{%
  \savebox{\mybox}{\hbox{\kern3pt$\scriptstyle#1$\kern3pt}}%
  \savebox{\mysim}{\hbox{$\sim$}}%
  \mathbin{\overset{#1}{\kern\z@\resizebox{\wd\mybox}{\ht\mysim}{$\sim$}}}%
}
\newcolumntype{C}[1]{>{\centering\let\newline\\\arraybackslash\hspace{0pt}}m{#1}}

\newcommand{\blind}{1}

\usepackage{ragged2e}

\newtheorem*{rem}{Remark}

\begin{document}

\def\spacingset#1{\renewcommand{\baselinestretch}%
{#1}\small\normalsize} \spacingset{1}


\if1\blind
{
  \centering{\bf\Large Controlling Sources of Inaccuracy in Stochastic Kriging}\\
  \vspace{0.2in}
  \centering{Wenjia Wang$^{a}$ and Benjamin Haaland$^{a,b}$\\
      \vspace{0.2in}
  \centering{$^a$Georgia Institute of Technology\\
    $^b$University of Utah}
} \fi

\if0\blind
{
  \bigskip
  \bigskip
  \bigskip
  \begin{center}
    {\bf\Large Controlling Sources of Inaccuracy in Stochastic Kriging}
\end{center}
  \medskip
} \fi

\begin{abstract}

Scientists and engineers commonly use simulation models to study real systems for which actual experimentation is costly, difficult, or impossible. Many simulations are stochastic in the sense that repeated runs with the same input configuration will result in different outputs. For expensive or time-consuming simulations, stochastic kriging \citep{ankenman} is commonly used to generate predictions for simulation model outputs subject to uncertainty due to both function approximation and stochastic variation. Here, we develop and justify a few guidelines for experimental design, which ensure accuracy of stochastic kriging emulators. We decompose error in stochastic kriging predictions into nominal, numeric, parameter estimation and parameter estimation numeric components and provide means to control each in terms of properties of the underlying experimental design. The design properties implied for each source of error are weakly conflicting and broad principles are proposed. In brief, space-filling properties ``small fill distance" and ``large separation distance" should balance with replication at distinct input configurations, with number of replications depending on the relative magnitudes of stochastic and process variability. Non-stationarity implies higher input density in more active regions, while regression functions imply a balance with traditional design properties. A few examples are presented to illustrate the results.
\end{abstract}

\justify

\noindent%
\vfill

\newpage
\spacingset{1.45} 

\section{Introduction}\label{intro}
In many situations actual physical experimentation is difficult or impossible, so scientists and engineers use  simulations, or \emph{computer experiments}, to study a system of interest.
For example,
\citet{mak} study a complex simulation model for turbulent flows in swirl injectors, which are used in a spectrum of propulsion and power-generation applications, under a range of geometric conditions,
\citet{burchell} estimate sexual transmissibility of human papillomavirus infection via a stochastic simulation model,
and \cite{moran} use the Cardiovascular Disease Policy Model to project cost-effectiveness of treating hypertension in the U.S. according to 2014 guidelines.
Commonly, these simulations require a cascade of complex calculations and simulator runs are \emph{expensive} relative to their information content.
To enable exploration of the relationship between inputs and outputs in the system of interest, a typical and apparently high-quality solution is to collect data at several input configurations, then build an inexpensive approximation, or \emph{emulator}, for the simulation. \par

In many cases, the data collected from the computer simulation is \emph{stochastic} in the sense that repeated runs with the same input configuration will have
different outputs, driven primarily by elements of the simulation model which are inherently stochastic.
Consider for example, the Coronary Heart Disease Policy Model, which is the simulation backbone underlying the cost-effectiveness study in \cite{moran}.
For each subject in a large cohort (the U.S. adult population), this model generates a simulated Markov trajectory through risk and event categories.
These trajectories involve, for each subject and time-increment, randomly assigning a new state according to a specified distribution.
Even if all the simulation settings, what we are calling inputs here, are unchanged, a new run of the simulation model will have slightly different random trajectories, and in turn slightly different outputs.
For emulation of stochastic computer experiments, the stochastic kriging model proposed in \cite{ankenman} has gained considerable traction as a quality approximation in a broad spectrum of real applications.
In the stochastic kriging model, output associated with each input is decomposed as the sum of a mean (Gaussian process) output and random (Gaussian) noise. \par


The accuracy of the stochastic kriging emulator depends strongly on how the data is collected \citep{Staum,haaland2011,haaland2014}.
Notably,
\cite{ankenman} provides a few useful results relating to mean squared prediction error (MSPE) \emph{integrated over the design space}
indicating that the distinct data sites should be relatively space-filling, while the number of replications is driven by the relative magnitudes of process and stochastic variability.
Unfortunately, these results are limited to stationary process covariance with no non-trivial regression functions in the process mean.
Further, no explicit consideration is given to very important experimental design impacts on numeric stability and parameter estimation (or numeric stability in parameter estimation).
A spectrum of practical sequential design heuristics for stochastic kriging are explored in \cite{chen}.


\cite{haaland2014} examine the qualitative features of high-quality experimental designs for building accurate Gaussian process emulators of \emph{deterministic} computer experiments.
For deterministic emulators, it is shown that the weakly conflicting space-filling properties ``small fill distance'' and ``large separation distance'' ensure well-controlled error.
Non-stationarity in the process's correlation decay indicates a higher density of input locations in regions with more quickly decaying correlation, while non-trivial regression functions indicate a balance between the space-filling properties and traditional design properties targeting small variances of least squares coefficient estimates.
In the common situation where correlation parameters are estimated within the Gaussian process framework, space-filling designs are slightly shifted to emphasize particular sizes and orientations of pairwise differences between input locations.


Here, we
seek to develop and justify overarching
principles of data collection for \emph{stochastic} kriging.
Importantly, the primary target here is a \emph{qualitative} indication of what type of designs might be expected to enable one to build an accurate model, not optimal design.
Throughout, we will call a design ``high-quality'' with respect to a particular component of error if the relevant error component is well-controlled, so that the error is small, but not necessarily optimal.
Broadly, the development here follows the framework and many of the results laid out in \cite{haaland2014}.
Throughout, results which extend in a relatively straightforward manner from the deterministic case to the stochastic case will be described in brief, at a high level with differences highlighted, while completely unique results and those for which extension is more complex will be described in more depth.


Inaccuracy in stochastic kriging will be decomposed into four components, nominal, numeric, parameter estimation, and parameter estimation numeric error.
The overall approach is to bound these four types of error in terms of experimental design properties.
It will be shown that the implied design characteristics for these four sources of error are \emph{weakly conflicting}.
In Section \ref{preliminaries}, the problem is formally stated, some notation provided,
and
several important well-known results stated. Then, in respective Sections \ref{nominal}, \ref{numeric}, \ref{parameterEst}, and \ref{parameterEstNumeric}, the nominal, numeric, parameter estimation, and parameter estimation numeric error are bounded.
Designs which are high-quality with respect to the provided bounds are discussed and a few examples are given, with consideration to  stationary and non-stationary cases as well as non-trivial regression functions.
In Section \ref{NumericEg}, a few numeric examples are presented, examining how the process and noise variability relate to the balance between space-fillingness and replication.
Conclusions and implications are discussed briefly in Section \ref{discussion}.

\section{Preliminaries}\label{preliminaries}
\subsection{Stochastic Kriging Model}\label{stochKrigingModel}
We consider the situation where a \emph{noisy} output $y(x)$ can be observed at an input configuration $x$ in a compact set $\Omega\subset\mathbb{R}^d$.
The output is noisy, or stochastic, in the sense that another run, or observation, at $x$ will give a different output value.
The noisy outputs are modeled as the sum of a deterministic function plus mean zero Gaussian noise.
That is,
\begin{align}\label{eq:a}
y(x)&=f(x)+\epsilon(x),
\end{align}
where
$\epsilon(x)\sim N(0,\sigma^2_{\tau_*}(x))$ and $\tau_*\in \mathbb{R}^{p_1}$ is a vector of parameters.
Throughout, we will annotate \emph{true} parameter values, which are not subject to estimation or any type of numeric error, with an $*$ whenever this distinction between the true parameter values and their estimated or noisy counterparts is useful. Notably, the noise components are taken as independent across both input locations and replications at the same input location, and we have suppressed the dependence of $y(x)$ and $\epsilon(x)$ on a random element, say $\upsilon$.
Further, since we expect nearby input locations to have similar noise variances to some degree in most practical situations,
we adopt a finite dimensional model for the noise variance function $\sigma^2_{\tau}(\cdot)$ throughout.
For example,
one might model the noise variance as a log-linear function of splines with coefficients $\tau$.
Following \cite{ankenman}, the deterministic component $f:\Omega \rightarrow \mathbb{R}$
is modeled as
a Gaussian process (GP) (see, for example, \cite{fang} and \cite{santner}), 
$f\sim {\rm GP}(h(\cdot)^T\beta_*, \Psi_{\theta_*}(\cdot,\cdot))$ for some fixed, known \emph{regression} functions $h:\Omega \rightarrow \mathbb{R}^q$ and a positive definite covariance function $\Psi_{\theta_*}(\cdot,\cdot)$. Here, the process mean and covariance depend on respective unknown parameters $\beta_*\in \mathbb{R}^q$ and $\theta_* \in \mathbb{R}^{p_2}$.
Let $\vartheta=(\beta^T, \theta^T, \tau^T)^T$ denote the vector consisting of all the parameters.
Throughout, the underlying mean function in the stochastic kriging model will be considered as the primary estimation target.

As shown in \cite{ankenman}, the best linear unbiased predictor (BLUP), as well as its MSPE, can be expressed for fixed variance-covariance parameters, $\theta$ and $\tau$, in terms of the \emph{distinct} data locations and the average output at each.
The likelihood of the unknown parameters given the data, on the other hand, depends
on the individual outputs,
or more concisely the first and second moments at each distinct input location,
as shown in \cite{Binois}.
Throughout, we will use
notation following \cite{Binois}. Let $\bar Y$ denote the vector of average responses at each of the $n$ distinct locations and $\bar X$ to denote the \emph{corresponding} distinct design locations. On the other hand, we will use $Y$ to denote the full vector of $m$ outputs (not averaged) and $X$ to denote the corresponding (potentially non-distinct) design locations.
For the $i^{\rm th}$ \emph{distinct} design location $x_i$, let $k_i$ denote the number of replications observed at $x_i$.
Then, the experimental design corresponding to the $i^{\rm th}$ component of $\bar Y$ can be described in terms of the pair $(x_i, k_i)$ for $i=1,...,n$, where $x_i\in \Omega$ denotes a distinct design point, and $k_i$ denotes the number of replicates at $x_i$. Let
\begin{align*}
\bar{y}(x_i)=\frac{1}{k_i}\sum_{j=1}^{k_i}y_j(x_i).
\end{align*}
 denote the sample mean at point $x_i$, where $y_j(x_i)$ denotes the $j^{\rm th}$ experiment at $x_i$, so that $\bar{Y}=(\bar{y}(x_1),\ldots,\bar{y}(x_n))^T$.
Similarly, let $\bar \Sigma_{\epsilon}={\rm diag}\{\sigma^2_\tau(x_1)/k_1,\ldots,\sigma^2_\tau(x_n)/k_n\}$ denote the diagonal matrix of marginal \emph{noise} variances of the components $\bar Y$, and let $\Sigma_{\epsilon}$ denote the diagonal matrix of marginal \emph{noise} variances of the components $Y$.

Throughout, we will assume that our underlying experimental design is a component of a sequence of experimental designs with a convergent large sample distribution, in the sense that the corresponding sequence of empirical cumulative distribution functions converge pointwise to a fixed cumulative distribution function.
This assumption will be used to approximately control several quantities in terms of the \emph{large sample} properties of our design.

If $\beta$ is unknown, but both $\theta$ and $\tau$ are known,
then the BLUP for $f$ at an arbitrary location of interest $x\in\Omega$ is \citep{Staum}
\begin{align}
\hat{f}_\vartheta(x)=h(x)^T\hat{\beta}+\Psi_{\theta}(x,\bar X)[\Psi_{\theta}(\bar X,\bar X)+\bar \Sigma_{\epsilon}]^{-1}(\bar{Y}-H(\bar X)\hat{\beta}),\label{nominalEmulator}
\end{align}
where 
$\hat{\beta}=(H(\bar X)^T[\Psi_{\theta}(\bar X,\bar X)+\bar \Sigma_{\epsilon}]^{-1}H(\bar X))^{-1}H(\bar X)^T[\Psi_{\theta}(\bar X,\bar X)+\bar \Sigma_{\epsilon}]^{-1}\bar{Y}$,
$H(\bar X)$ has $i^{\rm th}$ row $h(x_i)^T$ for distinct data location $x_i$, and $\Psi_{\theta}(A,B)$ has elements $\Psi_{\theta}(a_i,b_j)$.
Similarly, the BLUP (\ref{nominalEmulator}) has expected squared prediction error (conditional on the observed data), or mean squared prediction error,
\begin{align}
\Psi_{\theta}(x,x)-(h(x)^T, \Psi_{\theta}(x,\bar X))
\left(\begin{array}{cc}
0 & H(\bar X)^T\\
H(\bar X) & \Psi_{\theta}(\bar X,\bar X)+\bar \Sigma_{\epsilon}
\end{array}
\right)^{-1}
\left(\begin{array}{c}
h(x) \\
\Psi_{\theta}(\bar X,x)
\end{array}
\right).\label{eq:MSE}
\end{align}
Applying block matrix inverse results \citep{harville}, the MSPE (\ref{eq:MSE}) can be written as
\begin{gather}
\begin{split}
&\Psi_{\theta}(x,x)-\Psi_{\theta}(x,\bar X)[ \Psi_{\theta} (\bar X,\bar X )+\bar \Sigma_{\epsilon} ]^{-1}\Psi_{\theta}(\bar X,x)\\
              &+ (h(x)-H(\bar X)^T[ \Psi_{\theta} (\bar X,\bar X)+\bar \Sigma_{\epsilon} ]^{-1}\Psi_{\theta}(\bar X,x))^T\\
              &\quad\times(H(\bar X)^T[ \Psi_{\theta} (\bar X ,\bar X)+\bar \Sigma_{\epsilon} ]^{-1}H(\bar X))^{-1}\\
              &\quad\times (h(x)-H(\bar X)^T[ \Psi_{\theta} (\bar X,\bar X)+\bar \Sigma_{\epsilon} ]^{-1}\Psi_{\theta}(\bar X,x)). \label{eq:MSEBlock}
\end{split}
\end{gather}

At first glance, the stochastic kriging model, which assumes a Gaussian process mean with Gaussian noise, appears quite narrow and restrictive.
In fact, the model is not as restrictive as it appears.
In particular, if one believes that the target function $f$ lies in a reproducing kernel Hilbert space (say for example, $f$ has a fixed number of continuous partial derivatives), then a representer theorem \citep{scholkopf2001generalized} ensures that the solution to a very broad range of loss or likelihood-based penalized regression problems has the form given in (\ref{nominalEmulator}), although $\hat\beta$ would be estimated differently and the regularizing matrix $\bar\Sigma_\epsilon$ constructed differently,
depending on the loss or likelihood.
In practice, the stochastic kriging model is typically a high-accuracy non-parametric estimate of the underlying function $f$, and would represent a high-quality starting approximation for each of the three examples mentioned in the first paragraph of this article (turbulent flows, sexual transmissibility, and cardiovascular policy).

The stochastic kriging model, which is adapted to simulations with noisy outputs, differs from a kriging model for simulations with deterministic outputs only by the inclusion of $\bar \Sigma_{\epsilon}={\rm diag}\{\sigma^2_\tau(x_1)/k_1,\ldots,\sigma^2_\tau(x_n)/k_n\}$ in the BLUP and MSPE formulas above.
In a sense, the kriging model for deterministic simulations is a special case of the stochastic kriging model for which $\sigma^2_\tau(\cdot)\equiv 0$.
Many of the results developed below extend immediately to the kriging model for deterministic simulations by taking $\sigma^2_\tau(\cdot)\equiv 0$ or the number of replications at the $i^{\rm th}$ distinct input location $k_i\to\infty$ across $i$.
Similarly, many of the results developed in \cite{haaland2014} for deterministic kriging
translate directly to the stochastic kriging context  with only a cosmetic rework.
The aspects of parameter estimation error that relate to estimation of $\sigma^2_\tau(\cdot)$, of course, do not.

\subsection{Sources of Inaccuracy}\label{sources}
As stated in the final paragraph of Section \ref{intro}, inaccuracy in stochastic kriging will be decomposed into four components, nominal, numeric, parameter estimation, and parameter estimation numeric error.
The numeric emulator is, in a sense, the actual, tangible emulator, which is subject to parameter estimation error as well as numeric error in both emulator calculation and parameter estimation.
Let $\vartheta_*$, $\hat{\vartheta}$, and $\tilde{\vartheta}$ respectively denote the true parameters, estimated parameters \emph{not} subject to floating point errors, and estimated parameters subject to floating point error in both computation and optimization.
As noted previously, $*$ will be used throughout to annotate true parameter values.
Similarly, we will use $\hat{\cdot}$ and $\tilde{\cdot}$ to identify quantities subject to estimation and numeric error, respectively.
Similar to decompositions in \cite{haaland2011} and \cite{haaland2014}, the norm of the difference between the estimator of the unknown function and real function can be decomposed into nominal, numeric, parameter estimation, and parameter estimation numeric components using the triangle inequality as follows,
\begin{gather}
\begin{split}
\| f-\tilde{f}_{\tilde{\vartheta}} \| & =  \| f-\hat{f}_{\vartheta_{*}}+\hat{f}_{\vartheta_{*}}-\hat{f}_{\hat{\vartheta}}+\hat{f}_{\hat{\vartheta}} -\hat{f}_{\tilde{\vartheta}}+\hat{f}_{\tilde{\vartheta}}-\tilde{f}_{\tilde{\vartheta}}  \| \\
& \leqslant  \underbrace{\| f - \hat{f}_{\vartheta_{*}}\|}_{{\rm nominal}} + \underbrace{\| \hat{f}_{\vartheta_{*}}-\hat{f}_{\hat{\vartheta}} \|}_{{\rm parameter\; estimation}} + \underbrace{\| \hat{f}_{\hat{\vartheta}} -\hat{f}_{\tilde{\vartheta}} \|}_{{\rm parameter\; numeric}} + \underbrace{\| \hat{f}_{\tilde{\vartheta}}-\tilde{f}_{\tilde{\vartheta}}\|}_{{\rm numeric}}.\label{eq:decomposition}
\end{split}
\end{gather}
Here $\tilde{f}_{\vartheta}$ denotes the nominal emulator subject to floating point errors in calculation.
Nominal error refers to the difference between the target function $f$ and its \emph{idealized} approximation $\hat{f}_{\vartheta_{*}}$, which is not subject to floating point or parameter estimation error.
Numeric error refers to the difference between the computed emulator $\tilde{f}_{\tilde{\vartheta}}$, which is subject to floating point arithmetic, and an idealized version of the emulator which is not subject to floating point error in emulator computation $\hat{f}_{\tilde{\vartheta}}$.
Parameter estimation error represents the difference between emulators with the true
and estimated parameters, $\hat{f}_{\vartheta_{*}}$ and $\hat{f}_{\hat{\vartheta}}$, respectively.
Parameter estimation numeric error refers to the difference between the emulator with numerically estimated parameters under floating point arithmetic $\hat{f}_{\tilde{\vartheta}}$ and the emulator under an \emph{exactly} estimated parameter $\hat{f}_{\hat{\vartheta}}$.
While decomposition (\ref{eq:decomposition}) holds for any norm, here the $L_2(\Omega)$ norm will be the primary focus.
Taking the expectation (conditional on the data) of (\ref{eq:decomposition}) and applying Jensen's inequality and Fubini's theorem \citep{shao} gives
\begin{gather}
\begin{split}
\mathbb{E}\| f-\tilde{f}_{\tilde{\vartheta}} \| & \leqslant  \sqrt{ \int_{\Omega}\mathbb{E}(f(x) - \hat{f}_{\vartheta_{*}}(x))^2{\rm d}x} + \sqrt{ \int_{\Omega}\mathbb{E}(\hat{f}_{\vartheta_{*}}(x)-\hat{f}_{\hat{\vartheta}}(x))^2{\rm d}x}\\
&\quad+\sqrt{ \int_{\Omega}\mathbb{E}(\hat{f}_{\hat{\vartheta}}(x) -\hat{f}_{\tilde{\vartheta}}(x))^2{\rm d}x} + \sqrt{ \int_{\Omega}\mathbb{E}(\hat{f}_{\tilde{\vartheta}}(x)-\tilde{f}_{\tilde{\vartheta}}(x))^2{\rm d}x}.
\label{eq:decomposition2}
\end{split}
\end{gather}
Notice that the BLUP with parameter $\vartheta_*$ is the nominal emulator $\hat{f}_{\vartheta_{*}}$ in the first term in (\ref{eq:decomposition2}) above, while the
portion of the first term, bounding the nominal error above, under the square root and inside the integral, $\mathbb{E}(f(x) - \hat{f}_{\vartheta_{*}}(x))^2$, equals the MSPE (\ref{eq:MSE}).

\section{Nominal Error}\label{nominal}
For a particular design problem, we have two approaches to reduce MSPE. The first approach is to add more distinct input locations to reduce the distance between potential inputs and design points, the other is to take more experimental runs at a particular location to reduce the predictive variance at that location. Intuitively, if there is a cluster of design points, then the ${\rm MSPE}$ of the experimental design including the cluster is almost the same as the ${\rm MSPE}$ of the experimental design with multiple experiments at one of the points in this cluster.
Our intuition is correct, as a consequence of the continuity of matrix summation, inverses, and quadratic forms, as summarized in Proposition \ref{clusterProp} below.
\begin{proposition}\label{clusterProp}
Suppose $f\sim {\rm GP}(h(\cdot)^T\beta, \Psi_{\theta}(\cdot,\cdot))$, for some fixed, known functions $h(\cdot)$ and a positive definite function $\Psi_{\theta}(\cdot,\cdot)$, with stochastic observations generated by the stochastic kriging model described in Section \ref{stochKrigingModel}.
Let $X=(X_1,X_2)$, where $X_1=(x_1,x_2,\ldots,x_r)$ and
\begin{gather}
X'=(\underbrace{x^*,\ldots,x^*}_{r\;{\rm replications}},X_2). \nonumber
\end{gather}
If
$\sigma_\tau^2(\cdot)>0$ 
and $\sigma_\tau^2(\cdot)$, $h(\cdot)$, and $\Psi_{\theta}(\cdot,\cdot)$ are continuous, then MSPE$(x)\rightarrow$ MSPE$'(x)$ as $x_i\rightarrow x^*$ for $i=1,...,r$, where MSPE$(x)$ is the MSPE of the BLUP based on $X$ and MSPE$'(x)$ is the MSPE of the BLUP based on $X'$.
\end{proposition}

\begin{rem}
Proposition \ref{clusterProp} indicates that there is not much difference, in terms of the resulting MSPE's across the input space, between a design with evaluations at new points close to an existing design point and a design with additional evaluations at that same existing point.
In particular, if there are reasons other than strictly accuracy, such as simplicity or cost, to prefer a design with fewer distinct input locations, then such a design may perform almost as well as a design with more distinct input locations.
\end{rem}
%
%
%
A bound on the nominal error for the uppermost terms of the MSPE (\ref{eq:MSEBlock}), which provide the MSPE for a mean model with no regression functions, is provided in Theorem \ref{nominalThm}. A proof is given in Section \ref{nominalThmPrf} of Appendix.
Notably, the proof follows the strategy laid out in the proof of Theorem 3.1 in \cite{haaland2014}, with a few additional complexities in handling $\bar \Sigma_{\epsilon}$.
In fact, the deterministic kriging result in Theorem 3.1 of \cite{haaland2014}, can be obtained as a special case of the Theorem below by setting $\lambda_{{\rm max}}(\bar \Sigma_{\epsilon})=0$.
Throughout, we will use the notation $\lambda_{\rm max}(A)$ and $\lambda_{\rm min}(A)$ to denote the maximum and minimum eigenvalues of a symmetric matrix $A$.

\begin{theorem}\label{nominalThm}
Suppose $f\sim GP(0,\Psi_{\theta}(\cdot,\cdot))$ for a positive definite function $\Psi_{\theta}(\cdot,\cdot)$ with $\Psi_{\theta}(\cdot,\cdot)\geqslant 0$, stochastic observations are generated by the stochastic kriging model described in Section \ref{stochKrigingModel}, $(n-2) \sup_{u,v\in \Omega}\Psi_{\theta}(u,v) > \lambda_{{\rm max}}(\bar \Sigma_{\epsilon})$,
then the MSPE of $f$ has upper bound 
\begin{gather}
\begin{split}
             &\Psi_{\theta}(x,x)-2\Psi_{\theta}(x_i,x)+\Psi_{\theta}(x_i,x_i)-\frac{(\Psi_{\theta}(x_i,x)-\Psi_{\theta}(x_i,x_i))^2}{n \sup_{u,v\in \Omega}\Psi_{\theta}(u,v)+\lambda_{{\rm max}}(\bar \Sigma_{\epsilon})}\\
             & +\frac{\lambda_{{\rm max}}(\bar \Sigma_{\epsilon})(n \sup_{u,v\in \Omega}\Psi_{\theta}(u,v)+2(\Psi_{\theta}(x_i,x)-\Psi_{\theta}(x_i,x_i)))}{n \sup_{u,v\in \Omega}\Psi_{\theta}(u,v)+\lambda_{{\rm max}}(\bar \Sigma_{\epsilon})}.\label{eq:NominalAllResult}
\end{split}
\end{gather}
\end{theorem}

Two special cases are examined.
These cases respectively represent broadly applicable \emph{stationary} and \emph{non-stationary} covariance models for the \emph{process} $f$, and will be referred to as the Stationary Model and Non-Stationary Model.
These models provide a concrete structure within which we can gain a qualitative understanding of the design implications of both stationarity and non-stationarity.
In the upcoming development, the overall bound on the uppermost terms of (\ref{eq:MSEBlock}) will be expressed in terms of the maximum of local bounds,
\begin{align}
\sup_{x\in \Omega}\mathbb{E}\{f(x)-\hat{f}_{\vartheta}(x)\}^2=\max_i \sup_{x\in A_i}\mathbb{E}\{f(x)-\hat{f}_{\vartheta}(x)\}^2,\label{eq:NominalSupE}
\end{align}
where $\cup_i A_i=\Omega$. The maximum over $i$ in (\ref{eq:NominalSupE}) can be controlled by imposing a uniform bound over each of its components.
Below, $\varphi(\cdot)$ is a decreasing function of its non-negative argument and $\bar \Gamma$ is diagonal.

\subsection{Stationary Model}
Suppose
$\Psi_{\theta}(u,v)=\sigma^2\varphi(\|\Theta(u-v)\|_2)$ with $\bar \Sigma_{\epsilon}=\sigma^2\bar \Gamma$, where $\sigma\in \mathbb{R}_+$ is a parameter and $\Theta \in \mathbb{R}^{d \times d}$ is a non-singular matrix, which could be a parameter in its own right or a function of a lower dimensional parameter. Consider using the bound (\ref{eq:NominalAllResult}) as a guidepost for identifying the features of a high-quality experimental design.
Unlike in the deterministic case discussed in \cite{haaland2014}, in the stochastic kriging case, the denominator influences the bound  (\ref{eq:NominalAllResult}) through $\bar \Sigma_{\epsilon}$, inducing a balance between the variance at each point and the fill distance.
Notice that in (\ref{eq:NominalAllResult}), the bound is an increasing function of $\Psi_{\theta}(x,x)-\Psi_{\theta}(x_i,x)=\Psi_{\theta}(x_i,x_i)-\Psi_{\theta}(x_i,x)$.
Let $A_i=V_i(\Theta)$, the
Voronoi cell \citep{aurenhammer} anchored by distinct data point $x_i$, with respect to a Mahalanobis distance \citep{mahalanobis}
\begin{gather}
V_i(\Theta)=\{x\in\Omega:d_\Theta(x,x_i)\leqslant d_\Theta(x,x_j)\quad\forall j\ne i\},\nonumber
\end{gather}
where $d_\Theta(u,v)=\sqrt{(u-v)^T\Theta^T\Theta(u-v)}$.
On $A_i=V_i(\Theta)$,
the bound given by (\ref{eq:NominalAllResult})
can be bounded in terms of the smallest value of $\Psi_{\theta}(x_i,x)$,
which is attained for $x$ maximizing $d_{\Theta}(x_i,x)$.
Taking the maximum over $i$, and letting $\nu=\varphi(0)-\varphi(\max_i\sup_{x\in V_i(\Theta)}d_{\Theta}(x_i,x))$, (\ref{eq:NominalAllResult})
can be rewritten as
\begin{align}
& \sigma^2\bigg(2\nu-\frac{\nu^2}{n\varphi(0)+\lambda_{{\rm max}}(\bar \Gamma)}+\frac{\lambda_{{\rm max}}(\bar \Gamma)(n \varphi(0)+2\nu)}{n\varphi(0)+\lambda_{{\rm max}}(\bar \Gamma)}\bigg).\label{eq:NominalCase1Result}
\end{align}
Notice that
\begin{align*}
\max_i\sup_{x\in V_i(\Theta)}d_{\Theta}(x_i,x) = \sup_{x\in \Omega}\min_id_{\Theta}(x_i,x)
\end{align*}
is the fill distance with respect to the distance $d_{\Theta}$. Since (\ref{eq:NominalCase1Result}) is an increasing function of $\nu \in [0,\varphi(0)]$, the upper bound can be controlled by demanding the fill distance is small,
\emph{balanced with small largest element of} $\bar \Gamma$.


\begin{rem}
In the context of the stationary stochastic kriging model described above, experimental designs which balance small fill distance, with respect to the distance $d_\Theta$, for the distinct input locations with replication targeting uniformly small $\bar \Sigma_{\epsilon}$ ensure well-controlled nominal error.
\end{rem}


\subsection{Non-Stationary Model}\label{nonStatModel}
Here, we consider a relatively simple model of non-stationarity, adapted from \cite{ba2012composite}, which forms a good approximation in many practical situations.
In brief, the correlation decay is taken to be composed of more rapidly and more slowly decaying components, with the emphasis on the components depending on the input locations.
This allows the correlation to decay more quickly in some regions of the input space and more slowly in others.
This model of non-stationarity is reasonably well-suited to situations where the surface of interest is varying more quickly in some input regions and more slowly in others, and the model provides a structure for examining the design implications of this type of non-stationarity.


Suppose
$\Psi_{\Theta}(u,v)=\sigma^2(\omega_1(u)\omega_1(v)\varphi(\|\Theta_1(u-v)\|_2)+\omega_2(u)\omega_2(v)\varphi(\|\Theta_2(u-v)\|_2))$ with $\bar  \Sigma_{\epsilon}=\sigma^2\bar \Gamma$, where $\sigma \in \mathbb{R}_+$ is a parameter, and $\Theta_1, \Theta_2 \in \mathbb{R}^{d\times d}$ are non-singular matrices, either parameters in their own right or functions of lower dimensional parameters.
For the Non-Stationary Model case,  assume in addition $\omega_1(\cdot),\omega_2(\cdot)\ge 0$ have Lipschitz continuous derivatives on $\Omega$ with Lipschitz constants $k_1$ and $k_2$, respectively, $\omega_1^2(\cdot)+\omega_2^2(\cdot)=1$,
$\Theta_1, \Theta_2$ are non-singular, and $\lambda_{\rm max}(\Theta^T_1 \Xi^T_2 \Xi_2 \Theta_1)< 1$, where $\Xi_2=\Theta_2^{-1}$. The final assumption can be interpreted as $\varphi(\|\Theta_2(\cdot-\cdot)\|_2)$ is narrower than $\varphi(\|\Theta_1(\cdot-\cdot)\|_2)$.
For the Non-Stationary Model, we will localize the bounds over \emph{unions} of Voronoi cells
$V_i^*=V_i(\Theta_1)\cup V_i(\Theta_2)$.

Similar to the Stationary Model, we take the maximum over $i$, and let
\begin{align}
\nu=\varphi(0)-\inf_{x\in V^*_i}\{\omega_1(x)\omega_1(x_i)\varphi(\|\Theta_1(x-x_i)\|_2)+\omega_2(x)\omega_2(x_i)\varphi(\|\Theta_2(x-x_i)\|_2)\}.\label{eq:NominalNu2Result}
\end{align}
Then,
(\ref{eq:NominalAllResult}) again gives upper bound (\ref{eq:NominalCase1Result}).
Using Lipschitz continuity of $\omega_1(\cdot),\omega_2(\cdot)$ and Taylor's theorem \citep{nocedal2006numerical}, it can be shown that \citep{haaland2014},
\begin{align}
\nu & \leqslant \varphi(0) - (\omega_1^2(x_i)\varphi(\sup_{x\in V_i^*}d_{\Theta_1}(x_i,x)) + \omega_2^2(x_i)\varphi(\sup_{x\in V_i^*}d_{\Theta_1}(x_i,x))\}\nonumber\\
  & \qquad\qquad\quad- \varphi(0) (k_1+k_2)\max_i \sup_{x\in V_i^*}\|x-x_i\|_2).
\label{eq:NominalNu2Bound}
\end{align}
By plugging the right-hand side of (\ref{eq:NominalNu2Bound}) into (\ref{eq:NominalCase1Result}), we can obtain an upper bound, and corresponding guidepost for identifying features of a high-quality nominal error experimental design.

Following the development in \cite{haaland2014}, it can be shown that for fixed $\sigma_\tau^2(x_i)/k_i$, $i=1,\ldots,n$ (or equivalently $\bar\Gamma$),
(\ref{eq:NominalCase1Result}) is bounded uniformly over the design space by an experimental design with smaller \emph{union} of Voronoi cells, with respect to both $d_{\Theta_1}$ and $d_{\Theta_2}$, in regions with more emphasis on the quickly decaying correlation, and vice versa.
Similar to \cite{haaland2014}, the global and local correlation emphases are given concretely by
\begin{equation}
\begin{split}
&\omega_k(x_i)^2\left(\varphi\left(\sup_{x\in V^*_i}{ d}_{{\Theta}_1}(x_i,x)\right)-\varphi\left(\sup_{x\in V^*_i}{ d}_{{\Theta}_2}(x_i,x)\right)\right),\quad
k=1,2.\nonumber
\end{split}
\end{equation}

\begin{rem}
In the context of the non-stationary stochastic kriging model described above, experimental designs which balance smaller fill distances for the distinct input locations
in regions of the input space with more rapidly decaying correlation and larger fill distances in regions with more slowly decaying correlation, with replication targeting uniformly small $\bar \Sigma_{\epsilon}$ ensure well-controlled nominal error.
\end{rem}

\subsection{Regression Functions}
Next, we consider the lowermost terms in (\ref{eq:MSEBlock}), expressing the contribution of the regression terms to the overall accuracy. The regression terms can be bounded as
\begin{gather}\label{regFunupperbdNominal}
\begin{split}
&(h(x)-H(\bar X)^T[ \Psi_{\theta} (\bar X,\bar X)+\bar \Sigma_{\epsilon} ]^{-1}\Psi_{\theta}(\bar X,x))^T(H(\bar X)^T[ \Psi_{\theta} (\bar X,\bar X)+\bar \Sigma_{\epsilon} ]^{-1}H(\bar X))^{-1}\\
              & \quad\times (h(x)-H(\bar X)^T[ \Psi_{\theta} (\bar X,\bar X)+\bar \Sigma_{\epsilon} ]^{-1}\Psi_{\theta}(\bar X,x))\\
&\leqslant\sigma^2\frac{n\varphi(0)+\lambda_{\rm max}(\bar \Gamma)}{\lambda_{\rm min}(H(\bar X)^TH(\bar X))}\|h(x)-H(\bar X)^T[ \Psi_{\theta} (\bar X,\bar X)+\bar \Sigma_{\epsilon} ]^{-1}\Psi_{\theta}(\bar X,x)\|^2_2,
\end{split}
\end{gather}
by applying eigenvalue bounds on matrix quadratic forms, eigenvalue identities for matrix inverses, and Gershgorin's theorem \citep{varga} (see, for example, the development between equations (\ref{eq:NominalAll2}) and (\ref{gershgorin}) in the Appendix).

The term $\lambda_{\rm max}(\bar \Gamma)$ encourages balanced replication in the sense that it encourages a small maximum of $\sigma_\tau^2(x_i)/k_i$.
The term $\lambda_{\rm min}(H(\bar X)^TH(\bar X))$ in the denominator, on the other hand, encourages some degree of \emph{traditional} design properties. For example, linear regression functions would push input locations towards the edges or corners of the design space. On the other hand, the final term is the sum of squared errors for smoothed estimates of the regression functions and would be expected to be small in precisely the same situations when the topmost terms in (\ref{eq:MSEBlock}) are small, under the assumption that the regression functions can be well-approximated using the kernel $\Psi_\theta$ \citep{haaland2014}. That is, replication and traditional design properties need to be balanced with fill distance-based criteria.
Notably, the influence of non-trivial regression functions on what constitutes a high-quality design applies whether the process is stationary or non-stationary.
The corresponding design, adapted to a stationary or non-stationary covariance, is slightly pushed towards a traditional design targeting small variances of regression coefficients.

\begin{rem}
In the context of stochastic kriging models with non-trivial regression functions, experimental designs which balance space-filling properties, of the stationary or non-stationary variety as appropriate, replication targeting uniformly small $\bar\Sigma_\epsilon$, and traditional design properties targeting low-variance regression function coefficient estimates ensure well-controlled nominal error.
\end{rem}

\subsection{Example Designs}
Here, we seek to illustrate the type of designs indicated by the nominal error bounds, and provide a measure of corroboration for the qualitative features of good experimental designs that the bounds suggest.
For a given practical context and hypothetical values for the covariance parameters, the actual nominal error (\ref{eq:MSEBlock}) is computable, and could represent a component of a reasonable objective.

These designs are obtained by minimizing the nominal error bounds given by \eqref{eq:NominalCase1Result}, by plugging \eqref{eq:NominalNu2Bound} into \eqref{eq:NominalCase1Result}, and by taking the summation of \eqref{regFunupperbdNominal} and \eqref{eq:NominalCase1Result}, for the respective stationary, non-stationary, and non-trivial regression functions situations.
Since the noise variance is constant over the region, we need only consider the case where the number of replications at each distinct input location are equal.
The designs which minimize the upper bounds can then be obtained by minimizing over the number of replications.
In general, finding high-quality experimental designs is challenging, particularly when the objective function, here the relevant nominal error bound, is non-smooth and non-convex.
For a given number of replications, we can adopt the homotopy continuation \citep{eaves} procedure applied in \cite{haaland2014}.
In brief, we optimize the bounds over several iterations, slowly transitioning from an easier objective to the target objective.

Example high quality designs for stochastic kriging problems in the stationary situation, across a range of ratios $\sigma^2_\tau(x_i)/\sigma^2$,
and the non-stationary situation, as well as the stationary situation along with a constant and linear regression functions are shown in Figure \ref{nominalFigure}.
For the stationary cases shown in Panels 1-4, distinct design locations arrange themselves in a space-filling pattern, minimizing the fill distance. For the the non-stationary case shown in Panel 5, more distinct design locations are needed in portions of the input space with more emphasis on the more rapidly decaying correlation.
For the situation where a constant and linear regression functions are included with a stationary stochastic process variance, distinct design locations are pushed towards the corners of the input space, balancing space-filling and traditional design properties.
As the ratio of noise variance to functional variance $\sigma^2_\tau(x_i)/\sigma^2$ increases, more replications are needed at each distinct design location, moving from no replication when $\sigma^2_\tau(x_i)/\sigma^2=0.03$ to four replications when $\sigma^2_\tau(x_i)/\sigma^2=0.45$.


\begin{figure}[h!]
    \centering
        \includegraphics[height=4.5in]{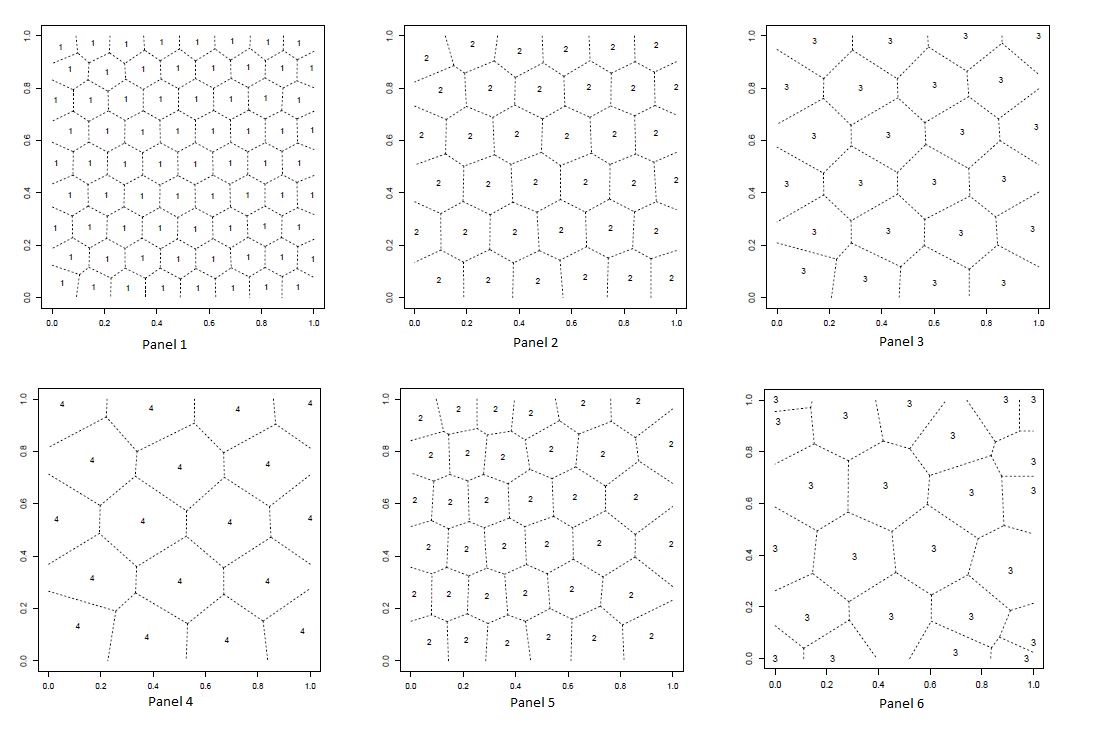}
    \caption{{\bf Panels 1-4:} Nominal error designs for \emph{stationary} correlation with $\varphi(d)=\exp\{-d^2\}$ and respective constant ratios of $\sigma_\tau^2(x_i)/\sigma^2$ equalling $0.03$, $0.10$, $0.25$, and $0.45$.
{\bf Panel 5:} Nominal error design for \emph{non-stationary} correlation with $\varphi(d)=\exp\{-d^2\}$, $\omega_1(x)=x_1$, $\Theta_1=I_2$, $\Theta_2=4I_2$, and ratio $\sigma_\tau^2(x_i)/\sigma^2$ of $0.10$.
{\bf Panel 6:} Nominal error design for \emph{stationary} correlation with $\varphi(d)=\exp\{-d^2\}$, and ratio $\sigma_\tau^2(x_i)/\sigma^2$ of $0.25$, along with a constant and two linear regression functions.
Design points annotated with number of replications throughout.}\label{nominalFigure}
\end{figure}

\section{Numeric Error}\label{numeric}
Numeric error comes from at least two sources. The first source is rounding error in the computer's representation of real numbers, and the second source is numeric solution to the parameter optimization problem.
In this section we develop bounds, in terms of properties of the experimental design, on the numeric error coming from the first numeric source of error, namely
$\| \hat{f}_{\tilde{\vartheta}}-\tilde{f}_{\tilde{\vartheta}}\|$. 
It can be shown that, similar to the non-stochastic kriging situation \citep{haaland2014}, increasing the number of data points always decreases the nominal error.
Unlike non-stochastic kriging, increasing the number of data points in the stochastic situation has far less ability to adversely affect numeric accuracy, particularly when $\sigma^2_\tau(x_i)$ is non-negligible.
It will be shown that the first source of numeric error can be controlled via the minimum eigenvalue
of $\Psi_{\theta}(\bar X,\bar X)+\bar \Sigma_{\epsilon}$, which has
\begin{align*}
\lambda_{\rm min}(\Psi_{\theta}(\bar X,\bar X)+\bar \Sigma_{\epsilon})&\geqslant\lambda_{\rm min}(\Psi_{\theta}(\bar X,\bar X))+\lambda_{\rm min}(\bar \Sigma_{\epsilon}).
\end{align*}

Numeric accuracy depends on the accuracy of floating point matrix manipulations. Commonly, computer and software have 15 digits of accuracy meaning roughly that
\begin{align*}
{\|\tilde{x}-x\|_2}/{\|x\|_2}\leqslant 10^{-15},
\end{align*}
where $x$ denotes the actual value and $\tilde{x}$ denotes the value that the computer stores.
Theorem \ref{numericThm} provides a bound on the numeric error of the stochastic GP emulator, which in turn requires relatively accurate
calculation of the functions $h$, $\bar Y$ and $\Psi$ (Assumption \ref{numericAssump}).

\begin{assumption}\label{numericAssump}
Assume $\kappa(\Psi_{\theta}(\bar X,\bar X)+\bar\Sigma_{\epsilon})=r/\delta$ with $r<1$, and
\begin{align*}
&\|h(x)-\tilde{h}(x)\|_2\leqslant \delta \|h(x)\|_2, \|\bar Y-\tilde{Y}\|_2\leqslant \delta \|\bar Y\|_2,\\
&\|\Psi_{\theta}(\bar X,\bar X)+\bar\Sigma_{\epsilon}-(\tilde{\Psi}_{\theta}(\bar X,\bar X)+\tilde{\Sigma}_{\epsilon})\|_2\leqslant \delta\|\Psi_{\theta}(\bar X,\bar X)+\bar\Sigma_{\epsilon}\|_2,\quad{\rm and}\\
& \|\Psi_{\theta}(x,\bar X)-\tilde{\Psi}_{\theta}(x,\bar X)\|_2\leqslant \delta\|\Psi_{\theta}(x,\bar X)\|_2,
\end{align*}
where $\delta>0$ quantifies the accuracy of computation in the computer.
Here, we use $\tilde{Y}$ and $\tilde{\Sigma}_{\epsilon}$ to denote the numeric representations of $\bar{Y}$ and $\bar{\Sigma}_{\epsilon}$, as opposed to the more cumbersome $\tilde{\bar{Y}}$ and $\tilde{\bar{\Sigma}}_{\epsilon}$.
\end{assumption}
The proof of Theorem \ref{numericThm} is essentially identical to the proof of Theorem 4.1 provided in the Appendix to \cite{haaland2014}, except with the additional $\bar{\Sigma}_{\epsilon}$ in the representation of the emulator (\ref{nominalEmulator}), so it is omitted for brevity.
\begin{theorem}\label{numericThm}
Suppose $f\sim {\rm GP}(h(\cdot)^T\beta, \Psi_{\theta}(\cdot,\cdot))$, for some fixed, known functions $h(\cdot)$ and a positive definite function $\Psi_{\theta}(\cdot,\cdot)$, with stochastic observations generated by the stochastic kriging model described in Section \ref{stochKrigingModel}.
For any fixed parameter estimate $\tilde{\vartheta}$, under Assumption \ref{numericAssump},
\begin{align*}
& \quad|\hat{f}_{\tilde{\vartheta}}(x)-\tilde{f}_{\tilde{\vartheta}}(x)|\nonumber\\
& \leqslant \delta\|h(x)\|_2\|\tilde{\beta}\|_2+\frac{2\delta}{1-r}\|\Psi_{\tilde{\theta}}(\bar  X,x)\|_2(\|H(\bar X)\|_2\|\tilde{\beta}\|_2+\|\bar Y\|_2)g(\Psi_{\theta}(\bar X,\bar X),\bar \Sigma_{\epsilon}),
\end{align*}
where 
\begin{gather}
\begin{split}
g(\Psi_{\theta}(\bar X,\bar X),\bar \Sigma_{\epsilon})&=\frac{1+\kappa(\Psi_{\theta}(\bar X,\bar X)+\bar \Sigma_{\epsilon})}{\lambda_{\min}(\Psi_{\theta}(\bar X,\bar X))+\lambda_{\min}(\bar \Sigma_{\epsilon})},\nonumber
\end{split}
\end{gather}
and $\kappa(\Psi_{\theta}(\bar X,\bar X)+\bar \Sigma_{\epsilon})$ denotes the condition number of $\Psi_{\theta}(\bar X,\bar X)+\bar \Sigma_{\epsilon}$.
\end{theorem}
Note that
\begin{gather}
\begin{split}
g(\Psi_{\theta}(\bar X,\bar X),\bar \Sigma_{\epsilon})&\leqslant\frac{1}{\lambda_{\min}(\Psi_{\theta}(\bar X,\bar X))+\lambda_{\min}(\bar \Sigma_{\epsilon})}\left(1+\frac{n \sup_{u,v\in \Omega}\Psi_{\theta}(u,v)+\lambda_{\max}(\bar \Sigma_{\epsilon})}{\lambda_{\min}(\Psi_{\theta}(\bar X,\bar X))+\lambda_{\min}(\bar \Sigma_{\epsilon})}\right),\label{eq:Numericg}
\end{split}
\end{gather}
where the inequality follows from Gershgorin's theorem \citep{varga}. See, 
equation (\ref{gershgorin}).


The norm $\|h(x)\|_2$ does not depend on the experimental design.
For experimental designs which are not too small and are high-quality with respect to \emph{parameter estimation numeric} properties, it will be shown in Section \ref{parameterEstNumeric}  that $\|\tilde{\beta}\|_2$ will approximately equal $\|\beta\|_2$.
Similarly, for experimental designs which are not too small and are high-quality with respect to \emph{nominal} properties, $\|\bar Y\|_2$ depends primarily on the sample size and large sample distribution of the inputs, as well as the target function $f$.
Further, for experimental designs which are not too small, the norms $\|\Psi_{\tilde{\theta}}(\bar X,x)\|_2$ and $\|H(\bar X)\|_2$
depend primarily on the sample size and large sample distribution of the inputs.
Thus, aside from $g(\Psi_{\theta}(\bar X,\bar X),\bar \Sigma_{\epsilon})$,
the other terms in the bound in the theorem influence the numeric error only weakly.
The bound depends on the experimental design primarily through $g(\Psi_{\theta}(\bar X,\bar X),\bar \Sigma_{\epsilon})$, which can be controlled via $\lambda_{\min}(\Psi_{\theta}(\bar X,\bar X))+\lambda_{\min}(\bar \Sigma_{\epsilon})$ as seen in (\ref{eq:Numericg}).
Unless $\lambda_{\min}(\bar \Sigma_{\epsilon})$ is \emph{very} near zero, the numeric error associated with generating predictions from a stochastic kriging model may be expected to be substantially less than in the deterministic case. On the other hand, when $\lambda_{\min}(\bar \Sigma_{\epsilon})$ \emph{is} very near zero, the numeric error in generating predictions would behave in a manner described in Section 3 of \cite{haaland2014}, favoring designs with \emph{well separated} distinct locations and, in the presence of non-stationarity, a greater (lesser) density of distinct input locations in sub-regions of the input space with more emphasis on local (global) correlation.
Within the framework described above, the relatively common practice in deterministic kriging of including a small so-called \emph{nugget} $\delta$, corresponding to $\bar \Sigma_{\epsilon}=\delta I_n$,
has the effect of greatly reducing numeric and parameter estimation numeric error, while (hopefully) only slightly increasing nominal and parameter estimation error.

\begin{rem}
In the context of the stochastic kriging model described above, numeric error is well-controlled for experimental designs with \emph{either} well-separated distinct input locations or $\lambda_{\min}(\bar \Sigma_{\epsilon})$ not too small.
\end{rem}

\section{Parameter Estimation Error}\label{parameterEst}
%
%

Throughout this section, the variance of the noise component $\sigma^2_{\tau}(x)$ is taken as
a continuously \emph{differentiable} function of
a finite dimensional, unknown parameter vector $\tau$, and maximum likelihood estimation is considered.
As noted in Section \ref{stochKrigingModel}, a finite dimensional parametrization of the noise variance $\sigma^2_{\tau}(x)$ can allow information sharing across nearby input locations, and accommodate small numbers of replications.
If, on the other hand, there are no extra assumptions on the noise variance (such as  similarity of nearby input locations), then the sample variance could be a sensible choice, as in \cite{ankenman}.
The likelihood of the parameters given the observed data can be expressed simply in terms of each individual output, and computed effiently via shortcut formulas provided in \cite{Binois}.
In this section, we will work with the full observation vector $Y$ and corresponding (potentially repeated) full design $X$.
Up to an additive constant, the log-likelihood is
\begin{align*}
l=-\frac{1}{2}\log\det[\Psi_{\theta}(X,X)+\Sigma_{\epsilon}]-\frac{1}{2}(Y-H(X)\beta)^T[\Psi_{\theta}(X,X)+\Sigma_{\epsilon}]^{-1}(Y-H(X)\beta).
\end{align*}
Let $\mathbb{E}$ denote the expectation conditional on $X$ and $Y$.
Then, for $n$ and $k_i$ not too small,
\begin{align}
\mathbb{E}\{\hat{f}_{\vartheta_{*}}(x)-\hat{f}_{\hat{\vartheta}}(x)\}^2 & \approx \frac{\partial \hat{f}(x)}{\partial \vartheta^T_*}\mathcal{I}(\vartheta_*)^{-1}\frac{\partial \hat{f}(x)}{\partial \vartheta_*}\label{pares},
\end{align}
where $\mathcal{I}(\vartheta_*)$ denotes the information matrix.
For approximation (\ref{pares}) to hold, we need the sequence of likelihood functions to become increasingly peaked.
For more details, see \cite{stein}.

In general, parameter estimates might be expected to affect the accuracy of Gaussian process regression models relatively weakly.
In fact, the order of approximation error will be the same across a wide range of parameter estimates, as long as the target function is in the reproducing kernel Hilbert space associated with the basic kernel \citep{haaland2011}.
The Fisher's information based error approximation in (\ref{pares}), while highly accurate only for large (and informative) samples, provides guidance for ensuring that the data we collect will enable construction of parameter estimates within this wide acceptable range.

%
%
%
%
For parameter estimation error, we have the following theorem, whose proof is provided in Section \ref{thm3prf} of Appendix.
Similar to Theorem \ref{nominalThm}, the proof of Theorem \ref{theopEstBound} follows the strategy laid out in the proof of Theorem 5.1 in \cite{haaland2014}, with a few additional complexities in handling the noise variance parameters $\tau$.
Once again, the deterministic kriging result in Theorem 5.1 of \cite{haaland2014}, can be obtained as a special case of the Theorem below by setting $\sigma^2_\tau(\cdot)=0$ and omitting the $c_4$ terms.
In the theorem, the Gaussian process covariance's parameters are separated as $\Psi_\theta(\cdot,\cdot)=\sigma^2\Phi_\rho(\cdot,\cdot)$.
\begin{theorem}\label{theopEstBound}
Let $f\sim {\rm GP}(h(\cdot)^T\beta, \sigma^2\Phi_\rho(\cdot,\cdot))$ for some fixed, known functions $h(\cdot)$ and positive definite function $\Phi_\rho(\cdot,\cdot)$, with stochastic observations generated by the stochastic kriging model described in Section \ref{stochKrigingModel}. Suppose $\hat{\vartheta}$ is the maximum likelihood estimator of the full set of unknown parameters $\vartheta=(\beta,\sigma^2, \rho, \tau)$.
Then, an approximate upper bound for $\mathbb{E}\{\hat{f}_{\vartheta_{*}}(x)-\hat{f}_{\hat{\vartheta}}(x)\}^2$ is given by
\begin{align}
& \frac{\sigma^2\|c_1\|^2_2(m \sup_{u,v\in \Omega}\Phi_{\rho}(u,v)+\lambda_{{\rm max}}(\Sigma_{\gamma}))}{ms_2}
                                          +\frac{\sigma^4\|c\|_2^2(m \sup_{u,v\in \Omega}\Phi_{\rho}(u,v)+\lambda_{{\rm max}}(\Sigma_{\gamma}))^2}{m^2s_1},\label{paramEstBound}
\end{align}
where
\begin{align*}
c_1  =&\frac{\partial \hat{f}(x)}{\partial \beta}
     = h(x)-H(X)^T(\Phi_\rho(X,X)+\Sigma_{\gamma})^{-1}\Phi_\rho(X,x),\nonumber\\
(c_3)_j  =&\frac{\partial \hat{f}(x)}{\partial \rho_j}=\bigg(\frac{\partial \Phi_\rho(x,X)}{\partial \rho_j}
                                               -\Phi_\rho(x,X)(\Phi_\rho(X,X)+\Sigma_{\gamma})^{-1}\frac{\partial \Phi_\rho(X,X)}{\partial \rho_j}\bigg)\nonumber\\
                                               &\quad\quad\quad\quad\times (\Phi_\rho(X,X)+\Sigma_{\gamma})^{-1}(Y-H(X){\beta}), \quad j =1,...,p_2,\nonumber\\
(c_4)_t  =&\frac{\partial \hat{f}(x)}{\partial \tau_t}
     =\Phi_\rho(x,X)(\Phi_\rho(X,X)+\Sigma_{\gamma})^{-1}{\rm diag}(\frac{\partial \gamma_1}{\partial \tau_t}I_{k_1},...,\frac{\partial \gamma_i}{\partial \tau_t}I_{k_i},...,\frac{\partial \gamma_n}{\partial \tau_t}I_{k_n})\\
                &\quad\quad\quad\quad\times(\Phi_\rho(X,X)+\Sigma_{\gamma})^{-1}(Y-H(X)\beta),\quad t=1,...,p_1,\nonumber\\
c  =&(c_3^T,c_4^T),
\end{align*}
where $(c_j)_i$ denotes the $i^{\rm th}$ element in vector $c_j$,
$\gamma_i=\sigma_\tau^2(x_i)/\sigma^2$, $\Sigma_{\gamma}={\rm diag}(\gamma_1I_{k_1},...,\gamma_nI_{k_n})$
$k_i$ is the number of replicates on $i^{\rm th}$ point, $m=\sum_{i=1}^{n}k_i$, and $s_1$ and $s_2$ are respectively defined in (\ref{defines1}) and (\ref{s2def}) in Appendix.
\end{theorem}


%

The upper bound is approximate in the sense that for a sequence of experimental designs with convergent maximum likelihood parameter estimates and convergent large sample distribution, the probability that the upper bound is violated by more than $\varepsilon>0$ goes to zero.
Notably, non-stochastic versions of this bound could be derived in terms of notions such as \emph{discrepancy} and \emph{total variation} or \emph{modulus of continuity} (see for example \cite{niederreiter}).

%
%

Following the development in
\cite{haaland2014}, both $\|c_1\|_2^2$ and $\|c_3\|_2^2$
involve \emph{interpolation} errors, for the regression functions and the derivatives of the Gaussian process covariance, respectively,
and these components would be expected to be
small for high quality nominal designs.
The remaining terms in $c_3$ are either well-controlled for high quality numeric designs, in the case of $(\Psi_{\theta}(X,X)+\Sigma_{\epsilon})^{-1}$, or depend only weakly on aspects of the experimental design beyond its size and large sample distribution, in the case of $Y-H(X)\beta$.
For $c_4$, we have the following proposition, whose proof is given in Section \ref{PfPropc4} of Appendix.
\begin{proposition}\label{Uppc4}
Under the conditions of Theorem \ref{theopEstBound},
\begin{align}\label{c4term}
|(c_4)_t| & \leqslant \frac{\|\Phi_\rho(x,\bar{X})\|_2\|\bar Y-H(\bar{X})\beta\|_2}{(\lambda_{{\rm min}}(\Phi_\rho(\bar{X},\bar{X})+\bar{\Sigma}_{\gamma}))^2}\max_{i:x_i\in\bar{X}}\bigg|\frac{1}{k_i}\frac{\partial \gamma_i}{\partial \tau_t}\bigg|.
\end{align}
\end{proposition}

The initial terms in (\ref{c4term}) are either well-controlled for high quality numeric designs, for $\lambda_{\rm min}(\Phi_\rho(\bar X,\bar X)+\bar \Sigma_{\gamma})$, or depend only weakly on aspects of the experimental design beyond its number of distinct locations and their large sample distribution, for $\|\Phi_\rho(x,\bar{X})\|_2$ and $\|\bar Y-H(\bar{X})\beta\|_2$.
The last term in (\ref{c4term}), $\max\left|\frac{1}{k_i}\frac{\partial \gamma_i}{\partial \tau_t}\right|$, encourages replication,
since it is a decreasing function of $k_i$.
Moreover, replication is more strongly encouraged near locations $x_i$ where $\gamma_i=\sigma_\tau^2(x_i)/\sigma^2$ is changing more rapidly with respect to one of the parameters $\tau_t$.
The term $s_2$ introduces a \emph{push} towards experimental design properties targeting reduction in variance of the regression function coefficients.

The term $s_1$ is somewhat more complex.
Let $W_1(x,y)=\Phi_\rho(x,y)+\sigma^2_\tau(x)/
\sigma^2\mathbb{I}_{\{x=y\}}$ and $\xi=(\rho,\;\tau)^T$. By (\ref{defines1}), $s_1 \geqslant 0$ and $s_1>0$ unless $\frac{\partial W_1(x,y)}{\partial\xi}a=W_1(x,y)b$ with probability 1 for some $(a,b)\neq 0$.
There are two parts to $\frac{\partial W_1(x,y)}{\partial\xi}$, $\frac{\partial \Phi_\rho(x,y)}{\partial\rho}$ and $\frac{\partial \sigma^2_\tau(x)}{\partial\tau}\mathbb{I}_{\{x=y\}}$.
Consider the \emph{distinct} and \emph{replicated} locations, $x\ne y$ and $x=y$,  separately.
The term $s_1$ will be large if two conditions are met.
First, the differences between distinct locations $\{x_i-x_j\}$ make $\frac{\partial \Phi_\rho(x_i,x_j)}{\partial \rho}$ far from zero, balanced with respect to a basis of $\mathbb{R}^{{\rm dim}\;\rho}$, and not collinear with $\Phi_\rho(x_i,x_j)$, similar to \cite{haaland2014}.
Second, the locations of replications make
$\frac{\partial \sigma^2_\tau(x_i)}{\partial\tau}$ far from zero, not collinear with $\Phi_\rho(x,x)+\sigma^2_\tau(x)/\sigma^2$, and \emph{balanced} in the sense that locations for which the derivative $\frac{\partial \sigma^2_\tau(x_i)}{\partial\tau}$ is \emph{small} in magnitude require \emph{more} replicates and \emph{vice versa}.
Notice that this encouragement of more replications where the derivative is smaller runs contrary to the influence of the term $\max\left|\frac{1}{k_i}\frac{\partial \gamma_i}{\partial \tau_t}\right|$ in $c_4$, which encourages more replications where the derivative is large in magnitude.
In our numeric examinations, the term encouraging more replication where the rate of change is large appears dominant.
Taken together, numeric studies suggest that the bound (\ref{paramEstBound}) is small for experimental designs whose distinct locations have good nominal and numeric properties, balanced with sufficient replications at each distinct data site, particularly at input locations where the noise variance is changing rapidly with respect one or more of the components of $\tau$.

\begin{rem}
In the context of the stochastic kriging model described above with parameters estimated via maximum likelihood, experimental designs with good nominal and numeric properties ensure well-controlled parameter estimation error.
\end{rem}




\subsection{Example Design}\label{paramEstExample}

Here again, we seek to illustrate the type of designs indicated by the parameter estimation error bounds, and provide a measure of corroboration for the qualitative features of good experimental designs that the bounds suggest.
For a given practical context and hypothetical values for the covariance parameters, the actual parameter estimation error (\ref{pares}) is approximately computable, and could represent a component of a reasonable objective.

Consider an example with $\Psi(d)=\exp(-d^Td)$, $\Psi_\rho(\cdot) = \Psi(\mbox{diag}\{\rho\}(\cdot))$, and $\rho = (1,1)^T$. In addition, suppose the stochastic error is given by $\sigma_\tau^2(x)=\tau\|x\|_2+0.04$, where $\tau$ is a parameter with true value $1$.
Suppose we want design points on $\Omega=[0,1]^2$. Since $\frac{\partial \sigma_\tau^2(x)}{\partial \tau} =\|x\|_2$, by (\ref{c4term}), a high quality experimental design should put more replicates on the locations that are far from zero. The total number of design points (may not be distinct locations) is 72, and the number of unique location is 24.
The corresponding design guided by the parameter estimation error bound \eqref{paramEstBound} is shown in Figure \ref{parEstFigure}.
Again, the locations of the distinct input setting are generated via the techniques described in \cite{haaland2014}, while the number of replications is taken so that the diagonal entries of $\bar{\Sigma}_\epsilon$ are nearly equal, to ensure well-controlled nominal error.
Notice that by balancing the stochastic error and $\frac{\partial \sigma_\tau^2(x)}{\partial \tau}$, the number of replicates are consistent with the contours of $\frac{\partial \sigma_\tau^2(x)}{\partial \tau}$, subject to edge effects.

\begin{figure}[h!]
    \centering
    \begin{subfigure}{.4\linewidth}
        \includegraphics[height=2in]{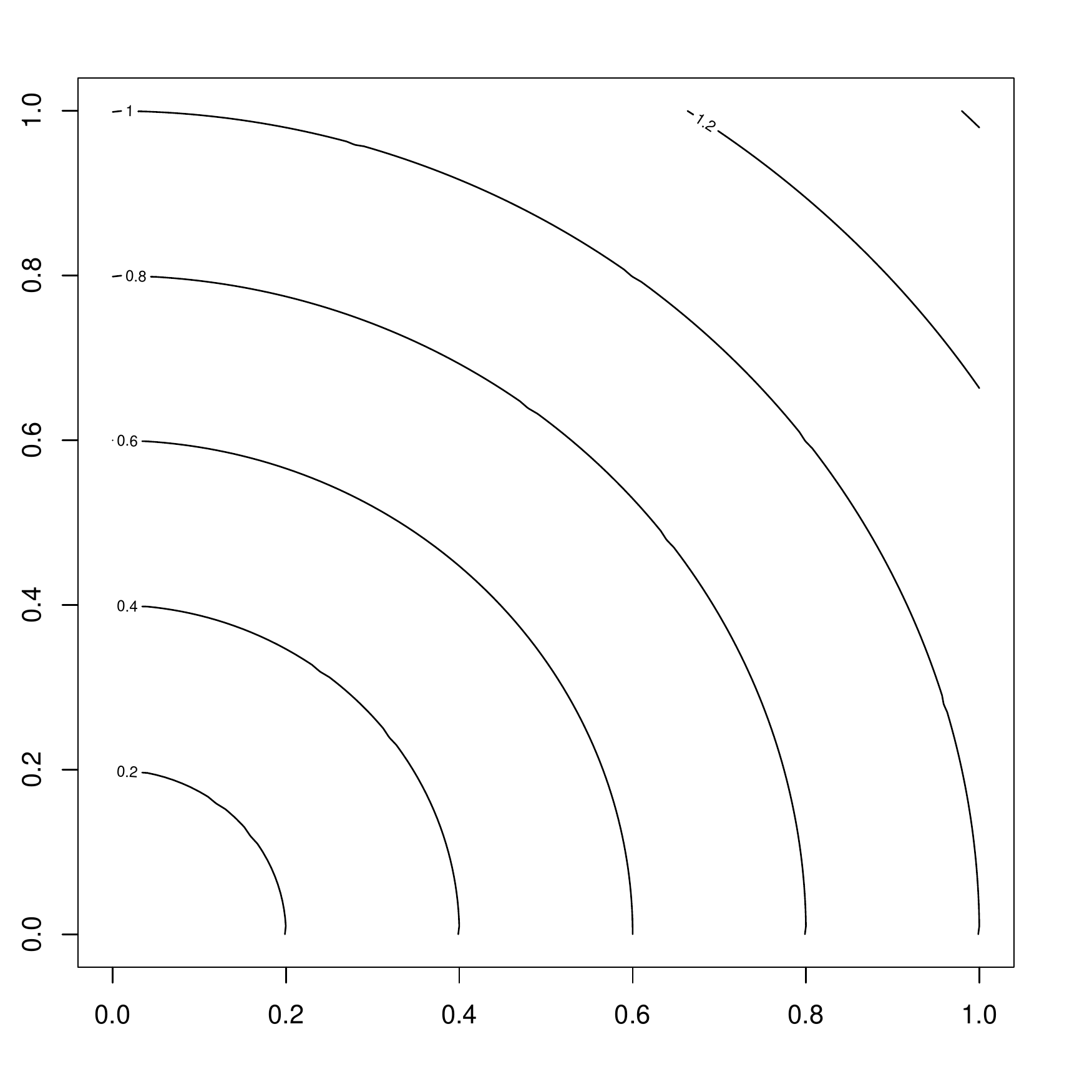}
    \end{subfigure}
    \begin{subfigure}{.4\linewidth}
        \includegraphics[height=2in]{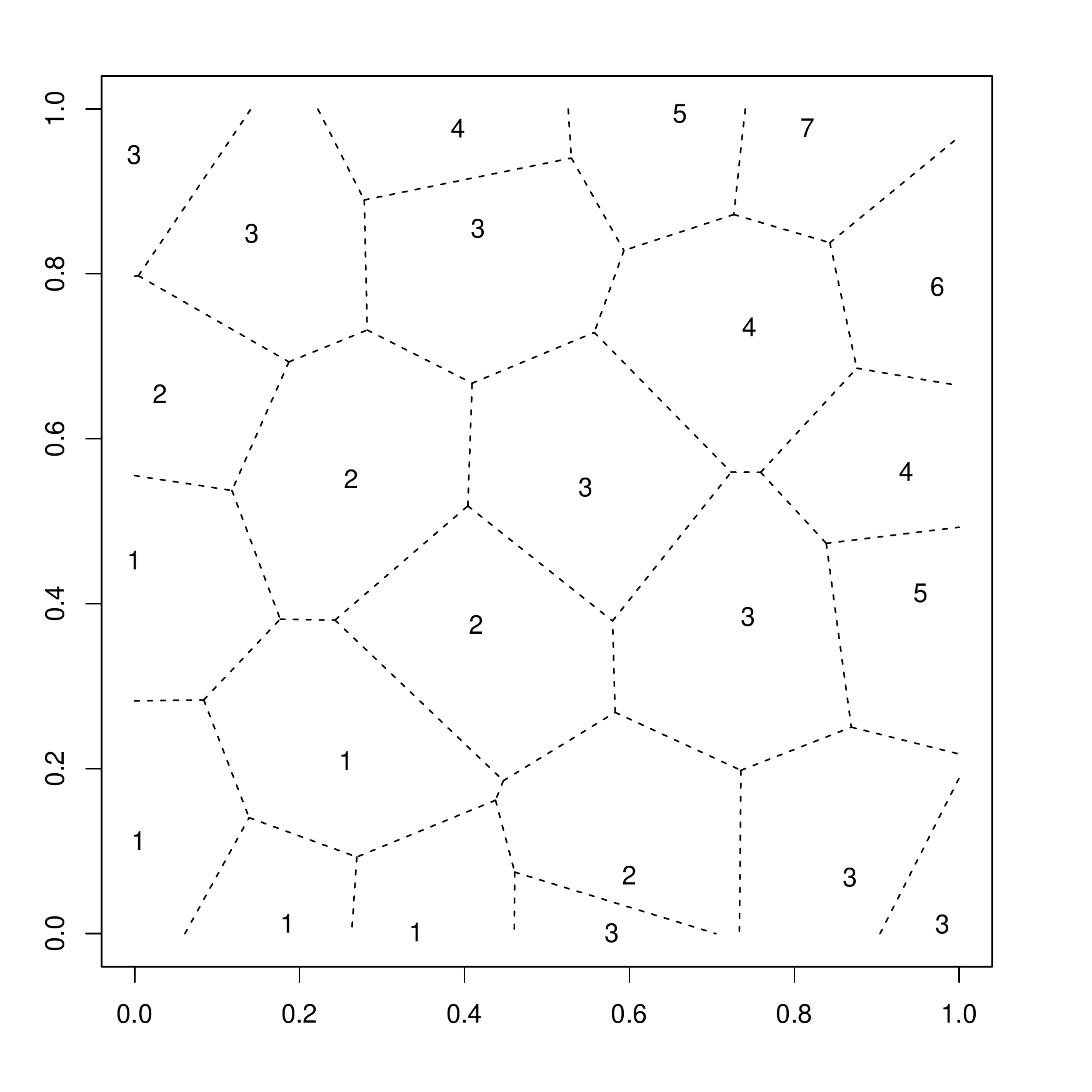}
    \end{subfigure}
    \caption{{\bf Left Panel:} Contour of $\frac{\partial \sigma_\tau^2(x)}{\partial \tau}$. {\bf Right Panel:} Parameter estimation error design with respect to the upper-bound in Theorem \ref{theopEstBound}. }\label{parEstFigure}
\end{figure}

\section{Parameter Estimation Numeric Error}\label{parameterEstNumeric}





In this section, the numeric error coming from numeric optimization of parameter estimates, the second source of numeric error, $\| \hat{f}_{\hat{\vartheta}} -\hat{f}_{\tilde{\vartheta}} \|_2$, is discussed.
Recall that $\hat{f}_\vartheta(x)=h(x)^T\beta+\Psi_{\theta}(x,\bar X)[\Psi_{\theta}(\bar X,\bar X)+\bar \Sigma_{\tau}]^{-1}(\bar y-H(\bar X)\beta)$, where each element of $\bar X$ denotes a distinct data location,
and $\bar\Sigma_{\tau}={\rm diag}(\sigma_\tau^2(x_1)/k_1,
\ldots,\sigma_\tau^2(x_n)/k_n)$.
Let $\tilde{A}=\Psi_{\tilde{\theta}}(\bar X,\bar X)+\bar \Sigma_{\tilde{\tau}}$ and $\hat{A}=\Psi_{\hat{\theta}}(\bar X,\bar X)+\bar \Sigma_{\hat{\tau}}$ denote the corresponding quantities subject to parameter estimation numeric error from numeric optimization and theoretical parameter estimates.
The below result links experimental design properties to parameter estimation numeric error. A proof is provided in Section  \ref{pfParaEstnum}
of Appendix.
\begin{theorem}\label{Thm:ParaEstNum}
Suppose $f\sim {\rm GP}(h(\cdot)^T\beta, \Psi_{\theta}(\cdot,\cdot))$
for some fixed, known function $h(\cdot)$ and positive definite function $\Psi_{\theta}(\cdot,\cdot)$, with stochastic observations generated by model (\ref{eq:a}).
Let $\tilde{\vartheta}$ denote the parameter we derive from numeric optimization and let $\hat{\vartheta}$ denote the true solution to the parameter optimization problem.
Let $\hat{f}_{\hat{\vartheta}}$ and $\hat{f}_{\tilde{\vartheta}}$ denote the BLUPs for $f$ with respective parameters ${\hat{\vartheta}}$ and ${\tilde{\vartheta}}$.
Then, under Assumptions \ref{numericAssump}, \ref{paramNumericAssump}, and \ref{assumpInParaNum},
\begin{align}\label{eq:ParaEstNucEqn}
|\hat{f}_{\hat{\vartheta}}(x)-\hat{f}_{\tilde{\vartheta}}(x)|
\leqslant & \frac{2\delta\kappa(\hat{A})}{(1-r)\lambda_{\min}(\hat{A})} \|\Psi_{\hat{\theta}}(\bar X,x)\|_2(\|f(\bar X)\|_2+\|H(\bar X)\|_2\|\hat{\beta}\|_2)\nonumber\\
& +2\delta\bigg(\kappa(\hat{A})\kappa(H(\bar X)^TH(\bar X))\bigg(1+(1+\delta)^2+\frac{(1+\delta)^2}{1-r}\kappa(\hat{A})\bigg)+1\bigg)\nonumber\\
& \times(\|h(x)\|_2+\frac{1+r}{(1-r)\lambda_{\min}(\hat{A})}\|H(\bar X)\|_2\|\Psi_{\hat{\theta}}(\bar X,x)\|_2)\|\hat{\beta}\|_2.
\end{align}
\end{theorem}
\begin{remark}
If Assumption \ref{assumpInParaNum} does not hold, we can still use Lemma \ref{lemmaLinearSystems} to derive an upper bound of $|\hat{f}_{\hat{\vartheta}}-\hat{f}_{\tilde{\vartheta}}|$, which is of order $\delta\kappa(\hat{A})^3$.
\end{remark}
Most of the terms above also appeared in Theorem \ref{numericThm}. The parameter estimation numeric error can also be controlled via $\lambda_{\min}(\Psi_{\theta}(\bar X,\bar X))+\lambda_{\min}(\bar \Sigma_{\epsilon})$ as seen in (\ref{eq:Numericg}).
See Section \ref{numeric} for a detailed discussion. The term $\kappa(H(\bar X)^TH(\bar X))$ requires some degree of \emph{traditional} design properties, as discussed in Section \ref{nominal}.
Consideration of parameter estimation numeric error further underscores the importance of well-conditionedness in the context of numeric linear algebra.

\begin{rem}
In the context of the stochastic kriging model described above with variance-covariance parameters estimated by numerically maximizing the likelihood, experimental designs with good numeric properties, slightly shifted towards good traditional design properties if non-trivial regression functions are included, ensure well-controlled parameter estimation numeric error.
\end{rem}


\section{Numeric Examples}\label{NumericEg}

In this section, we report simulation studies comparing designs with different numbers of replications.
Notably, we focus on the relationship between the number of replications at each distinct input location and the relative sizes of process and noise variation, potentially varying over the input space.
The relationship between the space-filling properties of the distinct input locations and emulator accuracy is examined empirically in \cite{haaland2014}.

It is worth noting again that the main objective here is to generate a qualitative description of which type of designs might be expected to perform well in the context of stochastic kriging.
For a particular application, the nominal and parameter estimation errors can be approximately computed (for a hypothesized parameter vector) and might form reasonable components of an overall design objective.
The numeric and parameter estimation numeric errors, on the other hand, cannot be computed, so the bounds given in Theorems \ref{numericThm}  and \ref{Thm:ParaEstNum} might form reasonable components of an overall design objective.


\subsection{Constant ratio of noise and process variance}
Take
$\Psi(u,v) = \exp(-\|u - v\|_2^2)$, $\sigma^2 = 1$, and space of interest $\Omega = [0,1]^2$. The total number of design points (potentially non-distinct) is set at $72$, and the number of replicates varied across $1, 2, 3$, and $4$, for 72, 36, 24, and 18 distinct locations.
Take $\epsilon(x) \sim N(0,\sigma_\epsilon^2)$ for all $x\in \Omega$.

For the initial study, set $\sigma^2_\epsilon$ to be 0.5, 0.1, and 0.01. Designs examined for the distinct input locations include the nominal designs shown in the first four panels of Figure \ref{nominalFigure}, $S$-optimal Latin hypercubes \citep{stocki2005method},
random Latin hypercubes \citep{stein1987large}, random uniform designs, and MaxPro designs \citep{joseph2015maximum}. First, 300 draws from the Gaussian process with mean zero and the correlation function $\Psi(\cdot, \cdot)$ are generated. For each draw, the observations based on the design and a 100 point random uniform testing set are made. Random errors drawn from $N(0,\sigma_\epsilon^2)$ are added to each of the observations. Based on the observations with random noise on the design points, predictions are generated on the testing set, and the maximum squared prediction error is computed. The \texttt{R} \citep{R} packages \texttt{lhs} \citep{carnell2016package} and \texttt{MaxPro} \citep{ba2015maxpro} were used for generating Latin hypercube and MaxPro designs. The average maximum squared prediction error over the 300 draws is calculated, and the results are reported in Table \ref{staEgTable}.

For a particular choice of experimental design strategy for the distinct input locations, we see an overall trend favoring replication as noise increases and space-fillingness as noise decreases.
Similar to \cite{haaland2014}, we see good performance for MaxPro designs \citep{joseph2015maximum} for the distinct input locations, as well as designs selected via the nominal error bound (\ref{eq:NominalCase1Result}).

\begin{table}
\centering

\begin{subtable}{0.7\linewidth}
\begin{tabular}{|c|c|c|c|c|}
\hline
\multicolumn{5}{|c|}{$\sigma_\epsilon^2 = 0.5$}\\
\hline
Design & rep = 4 & rep = 3 & rep = 2 & rep = 1  \\
\hline
nominal & 0.206 &  \textbf{0.202} &  0.221 &  0.212\\
optLHS & 0.236 & 0.229 &  \textbf{0.221} &  0.244\\
randLHS & 0.261 & 0.240 &  0.246 &  \textbf{0.217}\\
random & 0.295 & 0.278 &  0.249 &  \textbf{0.237}\\
MaxPro & \textbf{0.192} & 0.214  &  0.214 &  0.203\\

\hline
\multicolumn{5}{|c|}{$\sigma_\epsilon^2 = 0.1$}\\
\hline
Design & rep = 4 & rep = 3 & rep = 2 & rep = 1  \\
\hline
nominal & 0.071 &  0.071 &  0.071 &  \textbf{0.065}\\
optLHS & 0.088 & 0.083 &  0.079 &  \textbf{0.073}\\
randLHS & 0.109 & 0.092 &  0.091 &  \textbf{0.081}\\
random & 0.137 & 0.117 &  0.095 &  \textbf{0.084}\\
MaxPro & \textbf{0.059} & 0.063  &  0.067 &  0.067\\

\hline
\multicolumn{5}{|c|}{$\sigma_\epsilon^2 = 0.01$}\\
\hline
Design & rep = 4 & rep = 3 & rep = 2 & rep = 1  \\
\hline
nominal & 0.012 &  0.012 &  0.013 &  \textbf{0.012}\\
optLHS & 0.020 & 0.018 &  0.017 &  \textbf{0.014}\\
randLHS & 0.033 & 0.024 &  0.020 &  \textbf{0.017}\\
random & 0.047 & 0.036 &  0.025 &  \textbf{0.017}\\
MaxPro & 0.012 & 0.011  &  0.011 &  \textbf{0.010}\\
\hline
\end{tabular}
\end{subtable}
\caption{Average maximum squared prediction error for a spectrum of experimental designs across numbers of replications. The optimal average maximum squared prediction error (corresponding to a particular replication strategy) for each distinct input location space-filling design framework is in bold font.}
\label{staEgTable}
\end{table}



Next, we examine the
quality of the nominal error bound (\ref{eq:NominalCase1Result}), as well as any potential losses in accuracy due to following the guidance of the nominal error bounds in terms of the number of replications.
Here, the Gaussian process draws follow the same settings as the previous study. The designs examined here are optimal Latin hypercube and MaxPro. The total number of (potentially non-distinct) design points is set at 72 for all compared designs, with numbers of distinct locations in $\{6, 8, 9, 12, 18, 24, 36, 72\}$, and corresponding numbers of replicates in $\{12, 9, 8, 6, 4, 3, 2, 1\}$.
Noise standard deviations $\sigma_\epsilon$ are taken in $\{0.05, 0.35, 0.5\}$.
Comparisons of the nominal error bound (\ref{eq:NominalCase1Result}) to the average maximum squared prediction error over 300
draws of the Gaussian process are presented in Figure \ref{comparisonFigure}.
The figure illustrates the trade-off between replication and space-fillingness of the distinct design locations for a given sample size.
We see an initial trend of decreasing error as the number of distinct locations increases, followed by gradually increasing error as the number of replications becomes small.
The ideal balance of replication and space-fillingness of the distinct locations, favors replication for higher noise levels and favors space-fillingness for lower noise levels.
While this behavior is observed for both the bounds and actual error, the trend is quite marked for the bounds and more subtle for the actual error.



\begin{figure}[h!]
    \centering
    \begin{subfigure}{.23\linewidth}
        \includegraphics[trim=0 15 0 40, clip, height=1.4in]{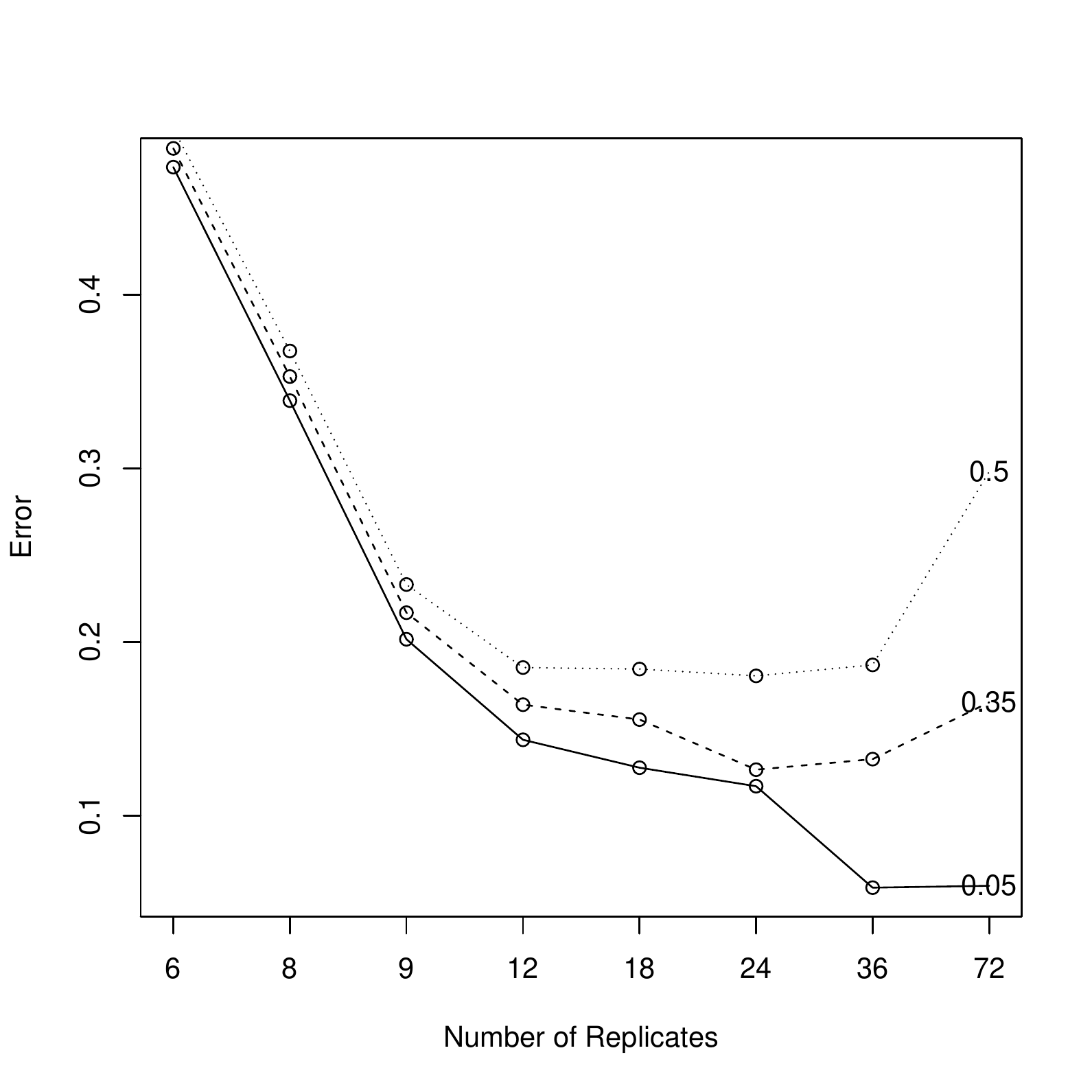}
    \end{subfigure}
    \begin{subfigure}{.23\linewidth}
        \includegraphics[trim=0 15 0 40, clip, height=1.4in]{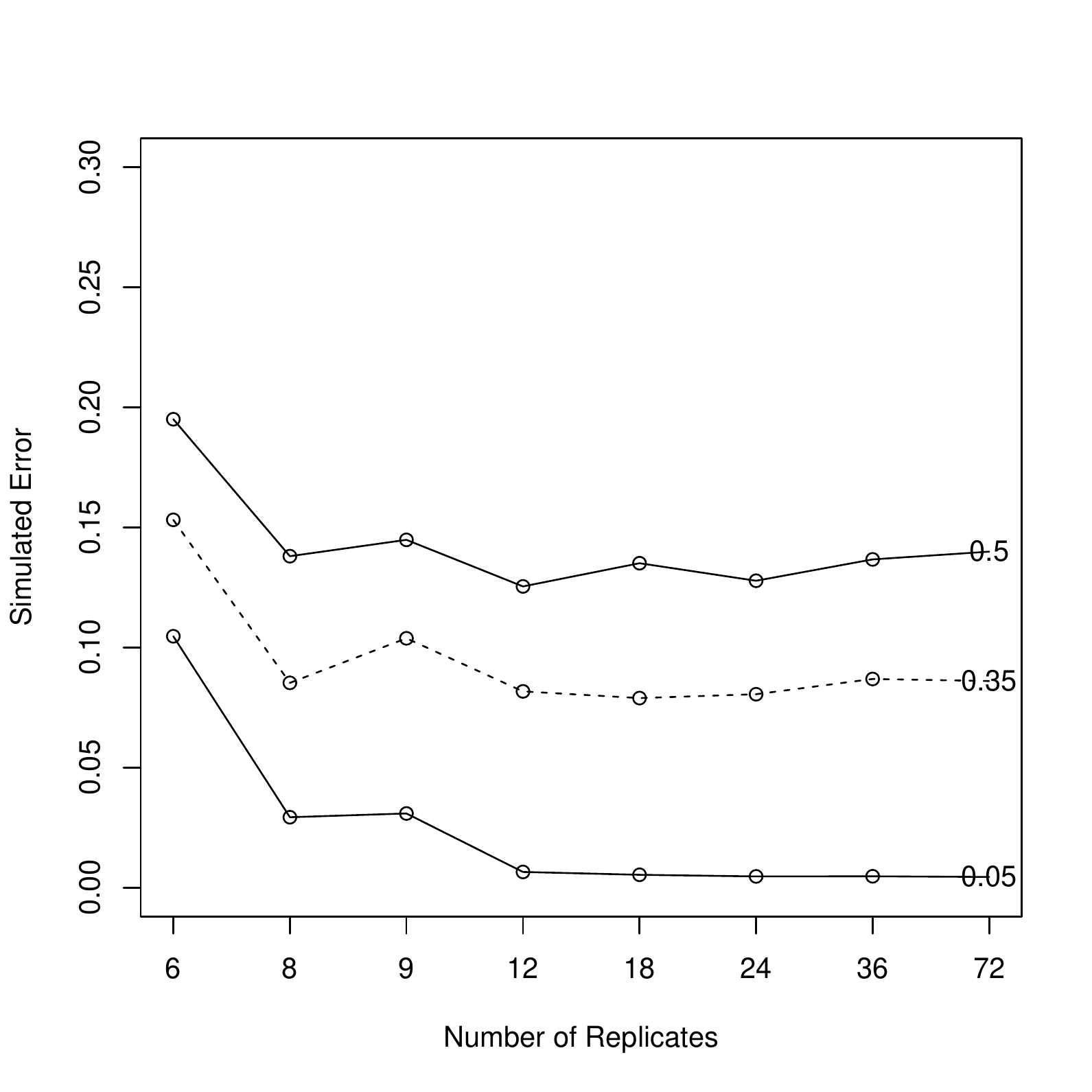}
    \end{subfigure}
    \begin{subfigure}{.23\linewidth}
        \includegraphics[trim=0 15 0 40, clip, height=1.4in]{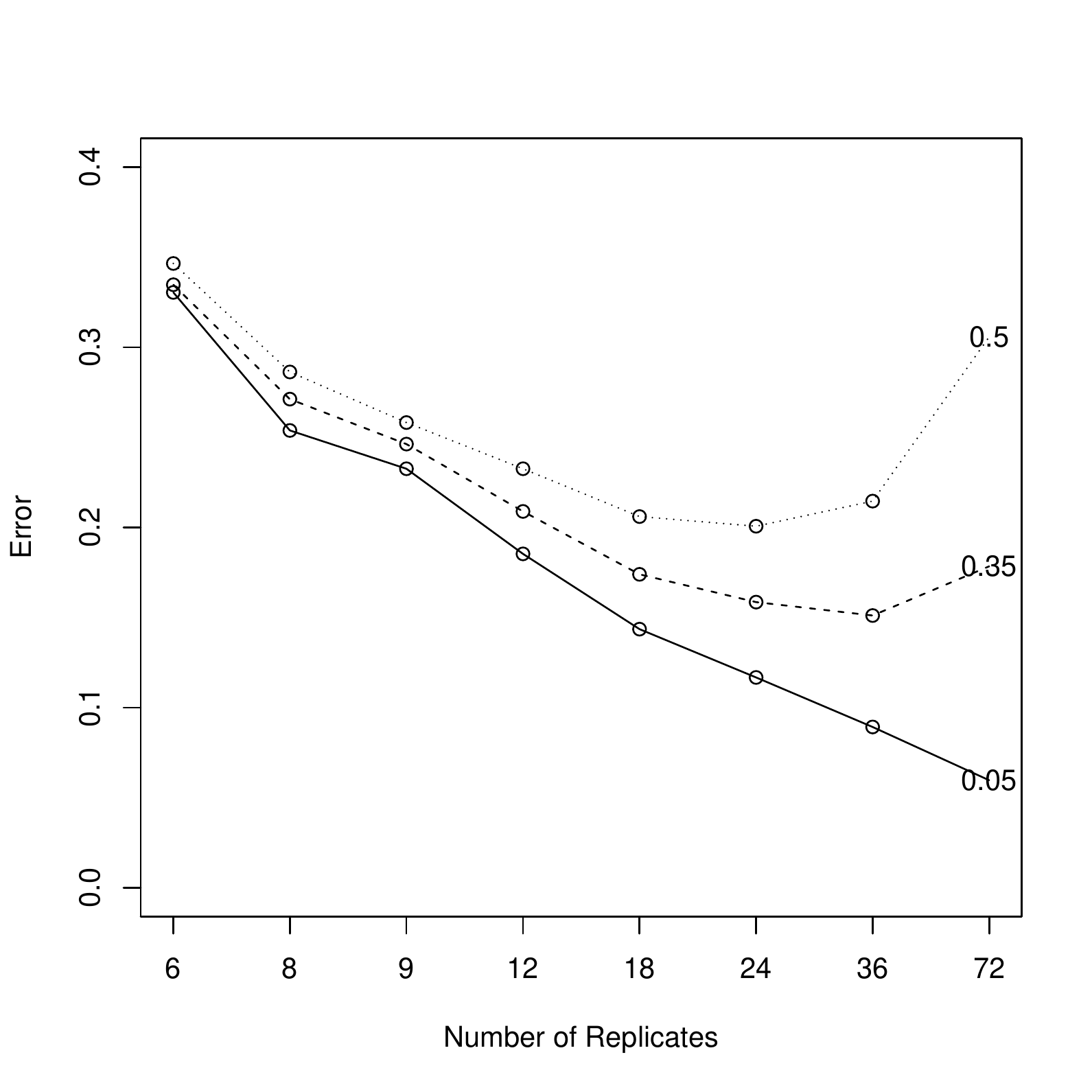}
    \end{subfigure}
    \begin{subfigure}{.23\linewidth}
        \includegraphics[trim=0 15 0 40, clip, height=1.4in]{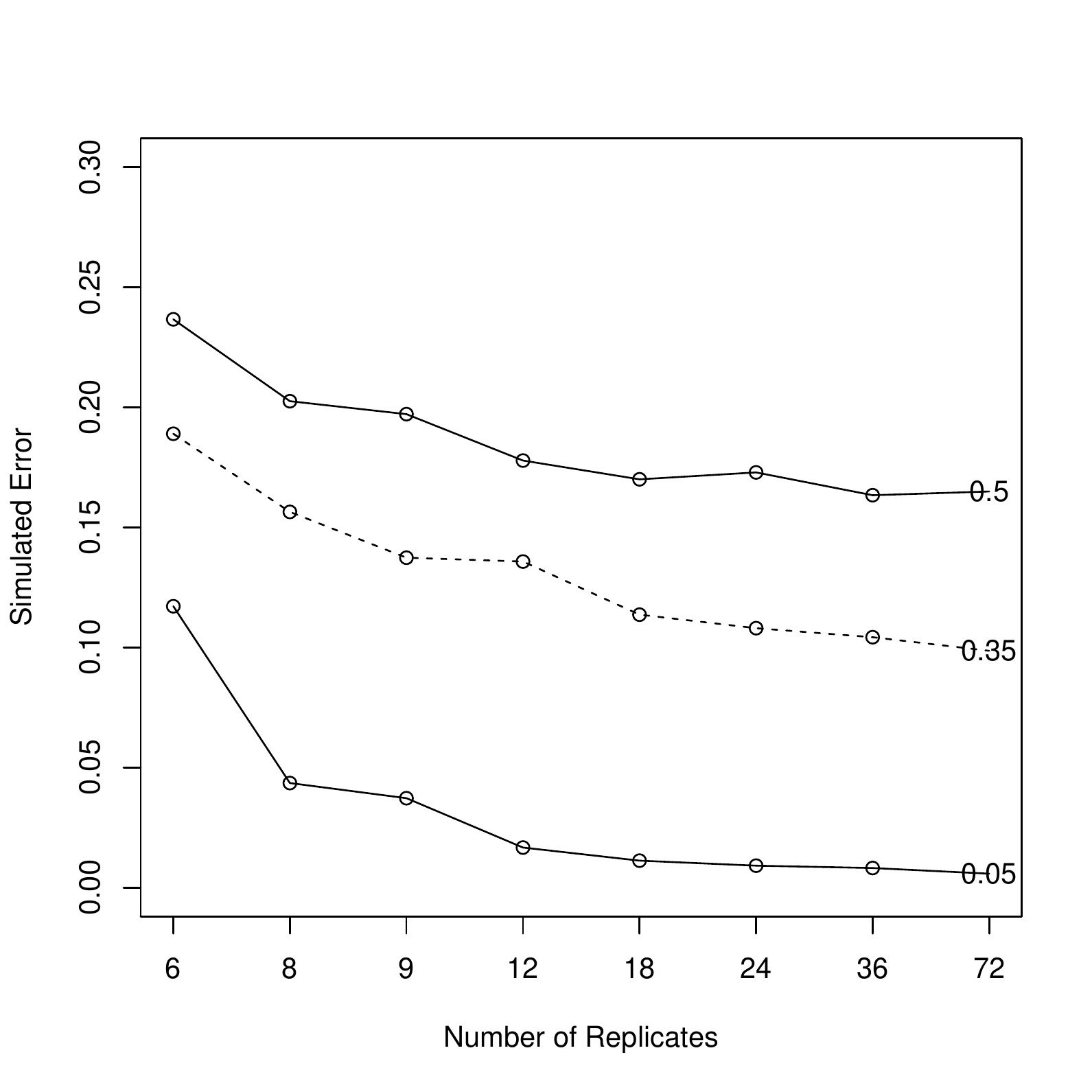}
    \end{subfigure}
    \caption{Comparisons of the nominal error bound \eqref{eq:NominalCase1Result} to the average maximum squared prediction error, where the total number of design points is 72. {\bf Left Panel:} Nominal error bound values for MaxPro designs. {\bf Middle Left Panel:} Average maximum squared prediction errors for MaxPro. {\bf Middle Right Panel:} Nominal error bound values for optimal Latin hypercube designs. {\bf Middle Left Panel:} Average maximum squared prediction errors for optimal Latin hypercube. }\label{comparisonFigure}
\end{figure}



Consider using the nominal error bound (\ref{eq:NominalCase1Result}) as guidance for choosing the number of replicates.
Here, we compare the average maximum prediction error under the best choice of replications to the average maximum prediction error under the number of replications suggested by the nominal bound. The noise standard deviations are taken to be $\sigma_\epsilon = 0.05k$ for $k = 1,\ldots,10$.
Relative and absolute differences in error are shown in Table \ref{tab:error}.
Results suggest that the bound provides useful guidance describing the qualities of a high-quality experimental design.

\begin{table}[h]
\centering
\begin{tabular}{|c|c|c|c|c|}
\multicolumn{1}{c}{}&\multicolumn{2}{c}{MaxPro}&\multicolumn{2}{c}{optLHS}\\
\hline
$\sigma_\epsilon$ & relative error & absolute error & relative error & absolute error   \\
\hline
0.05 & 0.132  & 0.00052& 0  & 0\\
0.10 & 0 & 0& 0  & 0\\
0.15 & 0.194 & 0.00401& 0.093 & 0.00285\\
0.20 & 0.049 & 0.00173& 0.040 & 0.00176\\
0.25 & 0.065 & 0.00331& 0 & 0\\
0.30 & 0.036 & 0.00209& 0.266 & 0.02271\\
0.35 & 0 & 0& 0.058 & 0.00577\\
0.40 & 0.003 & 0.00027& 0 & 0\\
0.45 & 0.065 & 0.00700& 0 & 0\\
0.50 & 0.076 & 0.00968& 0.030 & 0.00459\\
\hline
\end{tabular}
\caption{Relative and absolute loss of accuracy as compared to the optimal choice of replications, due to choosing the number of replications using the nominal error bound (\ref{eq:NominalCase1Result}), for MaxPro and optimal Latin hypercube designs.}
\label{tab:error}
\end{table}



\subsection{Input varying ratio of noise and process variance}
Next, we examine the model discussed in Section \ref{paramEstExample}, with noise level varying over the input space.
In particular, $\Psi(d)=\exp(-d^Td)$, $\Psi_\rho(\cdot) = \Psi(\mbox{diag}\{\rho\}(\cdot))$, and $\rho = (1,1)^T$, with stochastic error given by $\sigma_\tau^2(x)=\tau\|x\|_2+0.04$, where $\tau$ is a parameter with true value $1$.
Again, the input space is $\Omega=[0,1]^2$, and the total number of design points (potentially not distinct) is 72.
The parameter estimation design is provided by minimizing the upper bound provided in Theorem \ref{theopEstBound}.
Several designs for the distinct input locations including the nominal design provided in this paper, numeric designs  obtained along the lines described in \cite{haaland2014},
optimal Latin hypercubes \citep{stocki2005method}, random Latin hypercubes \citep{stein1987large}, random uniform designs, and MaxPro designs \citep{joseph2015maximum},
are considered.
We compare designs with equal replication at all distinct input locations and designs with unequal replications at the input locations as guided by \eqref{eq:NominalAllResult}, in which we require that the diagonal elements in $\bar \Sigma_\epsilon$ are nearly equal.
The number of unique locations used in the comparison are 18, 24, and 36. Then, we run 300 independent Gaussian process draws (with noise) and compare the average maximum squared prediction error of these processes. Results are shown in Table \ref{tab:varchange}.
In brief, accuracy is dramatically improved by varying the number of replications across the distinct input locations in the situation where the noise level varies across the input space.
This numeric study underscores the importance of appropriately handling input location dependent noise variance, as seen in the M/M/1 queue simulations in \cite{ankenman}, for both experimental design and modeling.

\begin{table}[h]
\centering
\begin{tabular}{|c|c|c|c|c|c|c|}
\multicolumn{1}{c}{}&\multicolumn{2}{c}{18 points} & \multicolumn{2}{c}{24 points} & \multicolumn{2}{c}{36 points}\\
\hline
Design & Varying & Const. & Varying & Const. & Varying & Const. \\
\hline
Nominal        & {\bf 0.139} & 0.182 & 0.146 & 0.192 & 0.164 & 0.212 \\
Numeric        & {\bf 0.108} & 0.155 & 0.128 & 0.185 & 0.159 & 0.193 \\
Parameter Est. & {\bf 0.113} & 0.141 & 0.125 & 0.170 & 0.125 & 0.190 \\
optLHS         & 0.157 & 0.209 & 0.144 & 0.209 & {\bf 0.129} & 0.198 \\
randLHS        & 0.176 & 0.210 & 0.162 & 0.228 & {\bf 0.148} & 0.225 \\
rand           & 0.229 & 0.267 & 0.184 & 0.239 & {\bf 0.160} & 0.227 \\
MaxPro         & {\bf 0.111} & 0.159 & 0.113 & 0.173 & 0.116 & 0.197 \\
\hline
\end{tabular}
\caption{Average maximum squared prediction error comparisons across number of distinct input locations and input varying replication vs. constant replication. The optimal average maximum squared prediction error (corresponding to a particular replication strategy) for each distinct input location space-filling design framework is in bold font.}
\label{tab:varchange}
\end{table}

\section{Discussion}\label{discussion}
We have
developed and justified guidelines for
ensuring accuracy of stochastic kriging predictors based on
experimental design.
By controlling nominal, numeric, parameter estimation and parameter estimation numeric sources of error, we
can control overall error in
stochastic kriging.
As in \cite{haaland2014}, the space-filling properties of the distinct design locations, ``small fill distance'' and ``large separation distance'', are largely non-conflicting with each of the sources of error.
Unlike \cite{haaland2014}, there is a trade-off between the number of replicates at each distinct design location and the space-filling properties of the distinct design locations.
This trade-off is reflected in the upper bounds for each of the four sources of errors.
The numeric error and parameter estimation numeric error are closely related to the condition number of $\Psi_{\theta}(\bar X,\bar X)+\bar \Sigma_{\epsilon}$, which becomes larger as more replicates or data locations are added.
Nominal and parameter estimation error, on the other hand, tend to encourage small fill distance.
Both the theoretical and numeric results suggest that the numeric and nominal properties of the distinct input locations, which are largely non-conflicting, as well as the number of replications at each, play a dominant role.

This work has several limitations.
Only upper bounds on the sources of error are considered.
There may be two designs with the same upper bound, where one is better than the other with respect to the expected error.
Our overall intention here is to provide a description of the qualitative features of high-quality experimental designs across a spectrum of common situations, rather than optimal design for a particular situation.
For a specific application, the nominal error (\ref{eq:MSEBlock}) and parameter estimation error (\ref{pares}) can actually be computed, for hypothesized parameter values, and could themselves represent components of a specific design objective.
The numeric error and parameter estimation numeric error, on the other hand, cannot be directly computed, so the respective bounds, provided in Theorems \ref{numericThm} and \ref{Thm:ParaEstNum}, might reasonably form the remaining components of the specific design objective.
We do not consider error from incorrectly using Gaussian process regression with maximum likelihood estimation to estimate the target function (model mis-specification).
From another perspective, the order of approximation error will be the same across a huge range of parameter estimates, as long as the target function is in the reproducing kernel Hilbert space associated with the basic kernel \citep{haaland2011}.
Projection design properties have not been explicitly discussed.
On the other hand, the results presented here indicate that if inert or inactive inputs are expected,  then the distinct design locations should be space-filling in lower-dimensional projections of the design.
Here, inert or inactive inputs refer to irrelevant input dimensions, over which the target function does not change.
Lastly, there are situations where a \emph{stochastic} emulator is need.
If the Gaussian \emph{noise} model fits the data well, then a stochastic emulator could be constructed by adding Gaussian noise with the estimated variance, $\sigma^2_{\hat{\tau}}(\cdot)$, where $\hat{\tau}$ denotes the maximum likelihood estimates of the noise variance parameters $\tau$.
Of course, the quality of this stochastic representation is directly linked to the quality of the noise parameter estimates, and the relevant design properties described in Section \ref{parameterEst} would be re-emphasized.
If the \emph{noise} model fits poorly, then perhaps a localized resampling of residuals could be useful.
The non-stationary model in Section \ref{nonStatModel} model would need to be altered slightly, for example $\Psi(x,y)=\sigma(x)\sigma(y)\varphi(\|x-y\|)$, to form a quality approximation for the setting where there is non-stationarity in the underlying process variability (as opposed to correlation).
Intuitively, we might expect this model of non-stationarity to necessitate more data in regions with more variability.
On another note, we might conceive of using these results in several ways including using the numeric and parameter estimation numeric bounds to remove points causing ill-conditioning.

In brief, these results provide further motivation and rationale for using one of several apparently high-quality, space-filling experimental designs for the distinct input locations, including but not limited to
\cite{iman1982distribution}, \cite{tang1993orthogonal}, \cite{owen1994controlling}, \cite{park1994optimal}, \cite{morris1995exploratory}, \cite{tang1998selecting}, \cite{ye1998orthogonal}, \cite{ye2000algorithmic}, \cite{jin2005efficient}, \cite{xu2011sudoku}, \cite{chen2014latin}, or \cite{joseph2015maximum},
particularly when there is no reason to expect non-stationarity in the process or noise.
While evidence of non-stationarity in process or noise variance, from an initial design perhaps, would indicate a varying density of distinct input locations or number of replications at distinct locations, respectively, precise characterization of this variation across the input space is challenging.
More generally, optimization of experimental designs is very challenging under many criteria, due to the high-dimensional and multi-modal nature of many of these problems.
For a situation with stationarity in both process and noise, a fixed number of replications across the design space, paired with one (or even a few) high-quality and computationally attractive space-filling designs for the distinct input locations, could be chosen in a computationally efficient manner and compared across numbers of replicates for a spectrum of plausible noise to process variance ratios.

\section*{Acknowledgements}
The authors would like to thank the Editor, Associate Editor and two anonymous Referees for their insightful comments and suggestions, which have greatly improved this manuscript.
The authors gratefully acknowledge funding from NSF DMS-1621722 and DMS-1739097.\\

\appendix
\section*{Appendix}
\renewcommand{\thesubsection}{\Alph{subsection}}
\renewcommand{\theequation}{\Alph{subsection}.\arabic{equation}}

\subsection{Proof of Theorem \ref{nominalThm}}\label{nominalThmPrf}
Consider a location of interest $x\in\Omega$ and the nearest \emph{design point} $x_i\in \bar X$. Let $\sigma_i^2 = \sigma^2_\tau(x_i)$
The uppermost terms in (\ref{eq:MSEBlock}) can be expressed as
\begin{gather}
\begin{split}
& \Psi_{\theta}(x,x)-\Psi_{\theta}(x,\bar X)[ \Psi_{\theta} (\bar X,\bar X)+\bar{\Sigma}_{\epsilon} ]^{-1}\Psi_{\theta}(\bar X,x)\\
              &=\Psi_{\theta}(x,x)\\
&\quad-[(\Psi_{\theta}(x,\bar X)-\Psi_{\theta}(x_i,\bar X)-\sigma_i^2e_i^T)[ \Psi_{\theta} (\bar X,\bar X)+\bar{\Sigma}_{\epsilon} ]^{-1}(\Psi_{\theta}(\bar X,x)-\Psi_{\theta}(\bar X,x_i)-\sigma_i^2e_i)\\
             &\quad+ 2(\Psi_{\theta}(x_i,\bar X)+\sigma_i^2e_i^T)[ \Psi_{\theta} (\bar X,\bar X)+\bar{\Sigma}_{\epsilon} ]^{-1}\Psi_{\theta}(\bar X,x)\\
             &\quad - (\Psi_{\theta}(x_i,\bar X)+\sigma_i^2e_i^T)[ \Psi_{\theta} (\bar X,\bar X)+\bar{\Sigma}_{\epsilon} ]^{-1}(\Psi_{\theta}(\bar X,x_i)+\sigma_i^2e_i)]\\
             &=\Psi_{\theta}(x,x)- 2e_i^T\Psi_{\theta}(\bar X,x)+e_i^T(\Psi_{\theta}(\bar X,x_i)+\sigma_i^2e_i)\\
             &\quad - (\Psi_{\theta}(x,\bar X)-\Psi_{\theta}(x_i,\bar X)-\sigma_i^2e_i^T)[ \Psi_{\theta} (\bar X,\bar X)+\bar{\Sigma}_{\epsilon} ]^{-1}(\Psi_{\theta}(\bar X,x)-\Psi_{\theta}(\bar X,x_i)-\sigma_i^2e_i)\\
&=\Psi_{\theta}(x,x)+\Psi_{\theta}(x_i,x_i)+\sigma_i^2-2\Psi_{\theta}(x_i,x)\\
&\quad - (\Psi_{\theta}(x,\bar X)-\Psi_{\theta}(x_i,\bar X)-\sigma_i^2e_i^T)[ \Psi_{\theta} (\bar X,\bar X)+\bar{\Sigma}_{\epsilon} ]^{-1}(\Psi_{\theta}(\bar X,x)-\Psi_{\theta}(\bar X,x_i)-\sigma_i^2e_i),
             \label{eq:NominalAll1}
\end{split}
\end{gather}
where $e_i$ denotes the $i^{\rm th}$ column of an $n\times n$ identity matrix.
The last term on the right-hand side of (\ref{eq:NominalAll1}) can be bounded as
\begin{align}
& -(\Psi_{\theta}(x,\bar X)-\Psi_{\theta}(x_i,\bar X)-\sigma_i^2e_i^T)[ \Psi_{\theta} (\bar X,\bar X)+\bar{\Sigma}_{\epsilon} ]^{-1}(\Psi_{\theta}(\bar X,x)-\Psi_{\theta}(\bar X,x_i)-\sigma_i^2e_i)\nonumber\\
             \leqslant & -\frac{\|\Psi_{\theta}(\bar X,x)-\Psi_{\theta}(\bar X,x_i)-\sigma_i^2e_i\|^2_2}{\lambda_{{\rm max}}[\Psi_{\theta} (\bar X,\bar X)+\bar{\Sigma}_{\epsilon}]}\nonumber\\
             \leqslant & -\frac{(\Psi_{\theta}(x_i,x)-\Psi_{\theta}(x_i,x_i)-\sigma_i^2)^2}{\lambda_{{\rm max}}[\Psi_{\theta} (\bar X,\bar X)+\bar{\Sigma}_{\epsilon}]}\nonumber\\
             \leqslant & -\frac{(\Psi_{\theta}(x_i,x)-\Psi_{\theta}(x_i,x_i))^2-2\sigma_i^2(\Psi_{\theta}(x_i,x)-\Psi_{\theta}(x_i,x_i))+\sigma_i^4}{\lambda_{{\rm max}}[\Psi_{\theta} (\bar X,\bar X)]+\lambda_{{\rm max}}(\bar{\Sigma}_{\epsilon})}\nonumber\\
             \leqslant
             & -\frac{(\Psi_{\theta}(x_i,x)-\Psi_{\theta}(x_i,x_i))^2-2\sigma_i^2(\Psi_{\theta}(x_i,x)-\Psi_{\theta}(x_i,x_i))+\sigma_i^4}{n \sup_{u,v\in \Omega}\Psi_{\theta}(u,v)+\lambda_{{\rm max}}(\bar{\Sigma}_{\epsilon})}\nonumber\\
             = & -\frac{(\Psi_{\theta}(x_i,x)-\Psi_{\theta}(x_i,x_i))^2}{n \sup_{u,v\in \Omega}\Psi_{\theta}(u,v)+\lambda_{{\rm max}}(\bar{\Sigma}_{\epsilon})}+\frac{2\sigma_i^2(\Psi_{\theta}(x_i,x)-\Psi_{\theta}(x_i,x_i))-\sigma_i^4}{n \sup_{u,v\in \Omega}\Psi_{\theta}(u,v)+\lambda_{{\rm max}}(\bar{\Sigma}_{\epsilon})},\label{eq:NominalAll2}
\end{align}
where the first inequality is true because for any vector $d$ and matrix $G$, $d^TG^{-1}d\geqslant \lambda_{\min}(G^{-1})\|d\|^2_2$ and $\lambda_{\min}(G^{-1})=1/\lambda_{\max}(G)$, the second inequality is true because the sum of squares $\|\cdot\|_2^2$ is larger than any one of its elements squared,
the third inequality is true because 
the maximum eigenvalue of a sum is at most the sum of the maximum eigenvalues,
and the final inequality is true because Gershgorin's theorem \citep{varga} implies
\begin{align}
\lambda_{\max}(\Psi_{\theta}(\bar X,\bar X))\leqslant \max_{j}\sum_{i=1}^n\Psi_{\theta}(x_i,x_j)\leqslant n \sup_{u,v\in \Omega}\Psi_{\theta}(u,v).\label{gershgorin}
\end{align}
Combining (\ref{eq:NominalAll1}) and (\ref{eq:NominalAll2}) gives
\begin{gather}
\begin{split}
& \Psi_{\theta}(x,x)-\Psi_{\theta}(x,\bar X)[ \Psi_{\theta} (\bar X,\bar X)+\bar{\Sigma}_{\epsilon} ]^{-1}\Psi_{\theta}(\bar X,x) \\
\leqslant & \Psi_{\theta}(x,x)+\Psi_{\theta}(x_i,x_i)-2\Psi_{\theta}(x_i,x) \\&-\frac{(\Psi_{\theta}(x_i,x) -\Psi_{\theta}(x_i,x_i))^2}{n \sup_{u,v\in \Omega}\Psi_{\theta}(u,v)+\lambda_{{\rm max}}(\bar{\Sigma}_{\epsilon})}+\sigma_i^2+\frac{2\sigma_i^2(\Psi_{\theta}(x_i,x)-\Psi_{\theta}(x_i,x_i))-\sigma_i^4}{n \sup_{u,v\in \Omega}\Psi_{\theta}(u,v)+\lambda_{{\rm max}}(\bar{\Sigma}_{\epsilon})}.\label{eq:NominalAllR1}
\end{split}
\end{gather}
Consider the concave, quadratic function
\begin{align*}
f_1(t) & = t+\frac{2t(\Psi_{\theta}(x_i,x)-\Psi_{\theta}(x_i,x_i))-t^2}{n \sup_{u,v\in \Omega}\Psi_{\theta}(u,v)+\lambda_{{\rm max}}(\bar{\Sigma}_{\epsilon})},
\end{align*}
where $t\in [0,\lambda_{{\rm max}}(\bar{\Sigma}_{\epsilon})]$. $f_1(\cdot)$ has axis of symmetry
\begin{align*}
t & = \frac{n \sup_{u,v\in \Omega}\Psi_{\theta}(u,v)+\lambda_{{\rm max}}(\bar{\Sigma}_{\epsilon})+2(\Psi_{\theta}(x_i,x)-\Psi_{\theta}(x_i,x_i))}{2}\nonumber\\
  & \geqslant \frac{(n-2) \sup_{u,v\in \Omega}\Psi_{\theta}(u,v)+\lambda_{{\rm max}}(\bar{\Sigma}_{\epsilon})}{2},
\end{align*}
where the last inequality is true because $\Psi_{\theta}(x_i,x)\geqslant 0 $ and $\Psi_{\theta}(x_i,x_i)<\sup_{u,v\in \Omega}\Psi_{\theta}(u,v)$.
If $(n-2) \sup_{u,v\in \Omega}\Psi_{\theta}(u,v) > \lambda_{{\rm max}}(\bar{\Sigma}_{\epsilon})$, then the axis of symmetry lies to the right of the interval $[0,\lambda_{{\rm max}}(\bar{\Sigma}_{\epsilon})]$ and $f_1(t)$ is increasing in $[0,\lambda_{{\rm max}}(\bar{\Sigma}_{\epsilon})]$. This indicates
\begin{gather}
\begin{split}
f_1(t) & \leqslant \lambda_{{\rm max}}(\bar{\Sigma}_{\epsilon})+\frac{2\lambda_{{\rm max}}(\bar{\Sigma}_{\epsilon})(\Psi_{\theta}(x_i,x)-\Psi_{\theta}(x_i,x_i))-\lambda_{{\rm max}}(\bar{\Sigma}_{\epsilon})^2}{n \sup_{u,v\in \Omega}\Psi_{\theta}(u,v)+\lambda_{{\rm max}}(\bar{\Sigma}_{\epsilon})}\\
& =\frac{\lambda_{{\rm max}}(\bar{\Sigma}_{\epsilon})(n \sup_{u,v\in \Omega}\Psi_{\theta}(u,v)+2(\Psi_{\theta}(x_i,x)-\Psi_{\theta}(x_i,x_i)))}{n \sup_{u,v\in \Omega}\Psi_{\theta}(u,v)+\lambda_{{\rm max}}(\bar{\Sigma}_{\epsilon})}.\label{eq:Nominalall3}
\end{split}
\end{gather}
Plugging (\ref{eq:Nominalall3}) into (\ref{eq:NominalAllR1}), gives the result.

\subsection{Proof of Theorem \ref{theopEstBound}}\label{thm3prf}
Here, the derivatives of the log-likelihood and emulator are expressed in terms of the equivalent parameters $\vartheta=(\beta^T,\sigma^2,\rho,\gamma)^T$, where $\gamma_i={\rm Var}(\epsilon(x_i))/\sigma^2$.
The vector of derivatives of the emulator with respect to the parameters $\frac{\partial \hat{f}(x)}{\partial \vartheta}$ has block components
\begin{align}
c_1=\frac{\partial \hat{f}(x)}{\partial \beta}&= h(x)-H(X)^T(\Phi_\rho(X,X)+\Sigma_{\gamma})^{-1}\Phi_\rho(X,x),\nonumber\\
c_2=\frac{\partial \hat{f}(x)}{\partial \sigma^2}&=0,\nonumber\\
(c_3)_j=\frac{\partial \hat{f}(x)}{\partial \rho_j}&=\frac{\partial \Phi_\rho(x,X)}{\partial \rho_j}(\Phi_\rho(X,X)+\Sigma_{\gamma})^{-1}(f(X)-H(X)\beta)\nonumber\\
                                               &\quad-\Phi_\rho(x,X)(\Phi_\rho(X,X)+\Sigma_{\gamma})^{-1}\frac{\partial \Phi_\rho(X,X)}{\partial \rho_j}(\Phi_\rho(X,X)+\Sigma_{\gamma})^{-1}(f(X)-H(X)\beta),\nonumber\\
(c_4)_t  =\frac{\partial \hat{f}(x)}{\partial \tau_t}
     &\Phi_\rho(x,X)(\Phi_\rho(X,X)+\Sigma_{\gamma})^{-1}{\rm diag}\left\{\frac{\partial \gamma_1}{\partial \tau_t}I_{k_1},...,\frac{\partial \gamma_i}{\partial \tau_t}I_{k_i},...,\frac{\partial \gamma_n}{\partial \tau_t}I_{k_n}\right\}\nonumber\\
                &\quad\times(\Phi_\rho(X,X)+\Sigma_{\gamma})^{-1}(y(X)-H(X)\beta),\nonumber
\end{align}
where $\Sigma_\gamma={\rm diag}(\gamma_1I_{k_1},...,\gamma_nI_{k_n})$.
The vector of derivatives of the log-likelihood with respect to the parameters $\frac{\partial l}{\partial \vartheta}$ has block components
\begin{align}
\frac{\partial l}{\partial \beta} & = \frac{1}{\sigma^2}(X)^T[\Phi_\rho(X,X)+\Sigma_{\gamma}]^{-1}(f(X)-H(X)\beta),\nonumber\\
\frac{\partial l}{\partial \sigma^2} & =-\frac{m}{2\sigma^2} +\frac{1}{2\sigma^4}(f(X)-H(X)\beta)^T(\Phi_\rho(X,X)+\Sigma_{\gamma})^{-1}(f(X)-H(X)\beta),\nonumber\\
\frac{\partial l}{\partial \rho_j} & = -\frac{1}{2}{\rm trace}\left([\Phi_\rho(X,X)+\Sigma_{\gamma}]^{-1} \frac{\partial \Phi_\rho(X,X)}{\partial \rho_j}\right)\nonumber\\
                                   &\quad + \frac{1}{2\sigma^2}(f(X)-H(X)\beta)^T[\Phi_\rho(X,X)+\Sigma_{\gamma}]^{-1}\frac{\partial \Phi_\rho(X,X)}{\partial \rho_j}[\Phi_\rho(X,X)+\Sigma_{\gamma}]^{-1}(f(X)-H(X)\beta),\nonumber\\
\frac{\partial l}{\partial \tau_t} & = -\frac{1}{2\sigma^2}{\rm trace}\left([\Phi_\rho(X,X)+\Sigma_{\gamma}]^{-1} {\rm diag}\left\{\frac{\partial \gamma_1}{\partial \tau_t}I_{k_1},...,\frac{\partial \gamma_i}{\partial \tau_t}I_{k_i},...,\frac{\partial \gamma_n}{\partial \tau_t}I_{k_n}\right\}\right)\nonumber\\
                                   &\quad + \frac{1}{2\sigma^4}(f(X)-H(X)\beta)^T[\Phi_\rho(X,X)+\Sigma_{\gamma}]^{-1}{\rm diag}(\frac{\partial \gamma_1}{\partial \tau_t}I_{k_1},...,\frac{\partial \gamma_i}{\partial \tau_i}I_{k_i},...,\frac{\partial \gamma_n}{\partial \tau_t}I_{k_n})\nonumber\\
                                   & \quad \times[\Phi_\rho(X,X)+\Sigma_{\gamma}]^{-1}(f(X)-H(X)\beta).\nonumber
\end{align}

So, the information matrix has block components
\begin{align*}
&\mathbb{E}-\frac{\partial^2 l}{\partial \beta^2}  = \frac{1}{\sigma^2}(X)^T[\Phi_\rho(X,X)+\Sigma_{\gamma}]^{-1}H(X),\\
&\mathbb{E}-\frac{\partial^2 l}{\partial \beta \partial \sigma^2}  = 0,\\
&\mathbb{E}-\frac{\partial^2 l}{\partial \beta \partial \tau_t}  = 0,\\
&\mathbb{E}-\frac{\partial^2 l}{\partial \beta \partial \rho_j}  = 0,\\
&\mathbb{E}-\frac{\partial^2 l}{\partial \sigma^2 \partial \sigma^2}  = \frac{m}{2\sigma^4}\\
&\mathbb{E}-\frac{\partial^2 l}{\partial \sigma^2 \partial \rho_j}  = \frac{1}{2}{\rm trace}\left([\Phi_\rho(X,X)+\Sigma_{\gamma}]^{-1}\frac{\partial \Phi_{\rho}(X,X)}{\partial
                                                           \rho_j}\right),\\
&\mathbb{E}-\frac{\partial^2 l}{\partial \rho_{j_1} \partial \rho_{j_2}}  = \frac{1}{2}{\rm trace}\left([\Phi_\rho(X,X)+\Sigma_{\gamma}]^{-1}\frac{\partial \Phi_\rho(X,X)}{\partial
                                                           \rho_{j_2}}[\Phi_\rho(X,X)+\Sigma_{\gamma}]^{-1}\frac{\partial \Phi_\rho(X,X)}{\partial
                                                           \rho_{j_1}}\right),\\
&\mathbb{E}-\frac{\partial^2 l}{\partial \tau_t \partial \sigma^2}  = \frac{1}{2\sigma^2}{\rm
                                                        trace}\left([\Phi_\rho(X,X)+\Sigma_{\gamma}]^{-1}{\rm diag}\left\{\frac{\partial \gamma_1}{\partial \tau_t}I_{k_1},...,\frac{\partial \gamma_i}{\partial \tau_t}I_{k_i},...,\frac{\partial \gamma_n}{\partial \tau_t}I_{k_n}\right\}\right),\\
&\mathbb{E}-\frac{\partial^2 l}{\partial \tau_t \partial \rho_j}  = \frac{1}{2}{\rm trace}\left([\Phi_\rho(X,X)+\Sigma_{\gamma}]^{-1}\frac{\partial
                                                        \Phi_\rho(X,X)}{\partial\rho_j}[\Phi_\rho(X,X)+\Sigma_{\gamma}]^{-1}\right.\\
&\quad\quad\quad\quad\quad\quad\quad\quad\quad\quad\left.\times{\rm diag}\left\{\frac{\partial \gamma_1}{\partial \tau_t}I_{k_1},...,\frac{\partial \gamma_i}{\partial \tau_t}I_{k_i},...,\frac{\partial \gamma_n}{\partial \tau_t}I_{k_n})\right\}\right),\nonumber\\
&\mathbb{E}-\frac{\partial^2 l}{\partial \tau_{t_1} \partial \tau_{t_2}}   = \frac{1}{2}{\rm trace}\left([\Phi_\rho(X,X)+\Sigma_{\gamma}]^{-1}{\rm diag}\left\{\frac{\partial \gamma_1}{\partial
                                                         \tau_{t_1}}I_{k_1},...,\frac{\partial \gamma_i}{\partial \tau_{t_1}}I_{k_i},...,\frac{\partial \gamma_n}{\partial \tau_{t_1}}I_{k_n}\right\}\right. & \\
& \quad\quad\quad\quad\quad\quad\quad\quad\quad\quad\left.\times [\Phi_\rho(X,X)+\Sigma_{\gamma}]^{-1}{\rm diag}\left\{\frac{\partial \gamma_1}{\partial \tau_{t_2}}I_{k_1},...,\frac{\partial \gamma_i}{\partial \tau_{t_2}}I_{k_i},...,\frac{\partial \gamma_n}{\partial \tau_{t_2}}I_{k_n}\right\}\right).
\end{align*}
Building from (\ref{pares}) and the block representations above gives
\begin{align*}
\mathbb{E}\{\hat{f}_{\vartheta_{*}}(x)-\hat{f}_{\hat{\vartheta}}(x)\}^2 & \approx (c_1^T,c_2^T,c_3^T,c_4^T)\left(\begin{array}{cc}a_{11}^{-1} & 0\\ 0 & \mathcal{I}^{-1}\\\end{array}\right)\left(\begin{array}{c}c_1\\c_2\\c_3\\c_4\\\end{array}\right)\\
& = c_1^T a_{11}^{-1} c_1+(c_2^T,c_3^T,c_4^T)\mathcal{I}^{-1}\left(\begin{array}{c}c_2\\c_3\\c_4\\\end{array}\right)\\
& = {\rm Part(I)}+{\rm Part(II)},
\end{align*}
where
\begin{align*}
a_{11}=\frac{\partial^2 l}{\partial \beta^2}\quad {\rm and}\quad
\mathcal{I}=\left(\begin{array}{cc}
\mathcal{I}_{11} & \mathcal{I}_{12}\\
\mathcal{I}_{21} & \mathcal{I}_{22}\\
\end{array}\right)
\end{align*}
with
\begin{align}
&\mathcal{I}_{11}=-\mathbb{E}\frac{\partial^2 l}{\partial \sigma^2 \partial \sigma^2},\quad
\mathcal{I}_{12}=-(\mathbb{E}\frac{\partial^2 l}{\partial \sigma^2 \partial \rho^T},\mathbb{E}\frac{\partial^2 l}{\partial \sigma^2\partial \tau^T}),\quad
\mathcal{I}_{21}=\mathcal{I}_{12}^T,\quad\mathcal{I}_{22}=\left(\begin{array}{cc}
D_{1} & D_{2}^T \\
D_{2}  & D_{3}\\
\end{array}\right),\nonumber\\
&D_{1}=-\mathbb{E}\frac{\partial^2 l}{\partial \rho \partial \rho^T}\quad
D_{2}=-\mathbb{E}\frac{\partial^2 l}{\partial \tau \partial \rho^T},\quad{\rm and}\quad
D_{3}=-\mathbb{E}\frac{\partial^2 l}{\partial \tau \partial \tau^T}.\label{eq:ParaEstD}
\end{align}
Applying block matrix inverse results \citep{harville} and noticing that $c_2=0$ gives
\begin{align*}
{\rm Part(II)}&=c^TB_1^{-1}c,
\end{align*}
where $B_1=\mathcal{I}_{22}-\mathcal{I}_{21}\mathcal{I}_{11}^{-1}\mathcal{I}_{21}$, and $c=(c_3^T,c_4^T)^T$.
With the aim of bounding ${\rm Part(II)}$, the following notation is introduced.
Let
\begin{align*}
a_j&={\rm vec}\left(\sigma^2\frac{\partial \Phi_\rho(X,X)}{\partial \rho_j}\right),\\
b_t&={\rm vec}\left({\rm diag}\left\{\frac{\partial \gamma_1}{\partial \tau_t}I_{k_1},...,\frac{\partial \gamma_i}{\partial \tau_t}I_{k_i},...,\frac{\partial \gamma_n}{\partial \tau_t}I_{k_n}\right\}\right),\nonumber\\
A_1&=(a_1,...,a_{|\rho|},b_1,...,b_{|\tau|}).
\end{align*}
Then,
\begin{align*}
B_1      & = \frac{1}{\sigma^4}A_1^T((\Phi_\rho(X,X)+\Sigma_{\gamma})^{-1}\otimes(\Phi_\rho(X,X)+\Sigma_{\gamma})^{-1}\nonumber\\
        &\quad-\frac{1}{n}{\rm vec}([\Phi_\rho(X,X)+\Sigma_{\gamma}]^{-1}){\rm vec}([\Phi_\rho(X,X)+\Sigma_{\gamma}]^{-1})^T)A_1.\nonumber
\end{align*}
For simplicity, let
\begin{align*}
W_1=\Phi_\rho(X,X)+\Sigma_{\gamma}\quad{\rm and}\quad
w=\frac{{\rm vec}(W_1)}{\|{\rm vec}(W_1)\|_2}.
\end{align*}
The matrix inside the quadratic form has eigenvector $w$ with corresponding eigenvalue 0.
Following the approach in \cite{haaland2014}, the minimum eigenvalue of $B_1$ can be bounded below by
\begin{align*}
\frac{1}{\sigma^4}\lambda_{\min}(A_1^T(I - ww^T)A_1) \times \lambda_2\bigg( (W_1^{-1}\otimes W_1^{-1}-\frac{1}{m}{\rm vec}(W_1^{-1}){\rm vec}(W_1^{-1})^T)\bigg),
\end{align*}
where $\lambda_2$ denotes the second smallest eigenvalue of its argument. Weyl's theorem \citep{ipsen} implies that the second smallest eigenvalue can be bounded below by
\begin{align*}
\lambda_{\min}(W_1^{-1}\otimes W_1^{-1}) = \frac{1}{\lambda_{\max}(\Psi_{\theta}(X,X)+\Sigma_{\epsilon})^2}\geqslant \frac{1}{( m \sup_{u,v\in \Omega}\Phi_\rho(u,v)+\lambda_{{\rm max}}(\Sigma_{\gamma}))^2}.
\end{align*}
For $\lambda_{\min}(A_1^T(I - ww^T)A_1)$, an approximate lower bound is given. Let $\xi = (\rho, \tau)$. Notice that
\begin{align}
&A_1^T(I - ww^T)A_1\nonumber\\
 & = \bigg[ \sum_{i,j} \frac{\partial W_1(x_i,x_j)}{\partial\xi}\frac{\partial W_1(x_i,x_j)}{\partial\xi^T}\nonumber\\
                     & \quad\quad- \frac{1}{\|{\rm vec}(W_1)\|_2^2}\bigg(\sum_{i,j} \frac{\partial W_1(x_i,x_j)}{\partial\xi}W_1(x_i,x_j)\bigg)\bigg(\sum_{i,j} \frac{\partial W_1(x_i,x_j)}{\partial\xi^T}W_1(x_i,x_j)\bigg)\bigg]\nonumber\\
                     & \approx m^2\bigg[ \int \frac{\partial W_1(x,y)}{\partial\xi}\frac{\partial W_1(x,y)}{\partial\xi^T}dF^2(x,y)\nonumber\\
                     & \quad\quad\quad- \frac{1}{\|W_1\|_{L_2(F^2)}^2}\bigg(\int \frac{\partial W_1(x,y)}{\partial\xi}W_1(x,y)dF^2(x,y)\bigg)\bigg(\int \frac{\partial W_1(y)}{\partial\xi^T}W_1(x,y)dF^2(x,y)\bigg)\bigg]\nonumber\\
                     & \succeq m^2 s_1, \label{defines1}
\end{align}
where $W_1(x,y)=\Phi_\rho(x,y)+\frac{\sigma^2_\tau(x)}{\sigma^2}\mathbb{I}_{\{x=y\}}$ and $F^2$ denotes the large sample distribution of point pairs. Applying a version of the Cauchy-Schwarz inequality for random vectors \citep{tripathi}, gives $s_1 \geqslant 0$ with $s_1>0$ unless
\begin{align*}
\frac{\partial W_1(x,y)}{\partial\xi}a=W_1(x,y)b
\end{align*} with probability 1 with respect to large sample distribution of point pairs $F^2$ for some vectors $a$ and $b$. So, ${\rm Part(II)}$ has approximate upper bound
\begin{align}
\frac{\sigma^4(m \sup_{u,v\in \Omega}\Phi_\rho(u,v)+\lambda_{{\rm max}}(\Sigma_{\gamma}))^2}{m^2s_1}\|c\|_2^2.\label{eq:AppendixPartII}
\end{align}
Also,
\begin{align}
{\rm Part(I)} \leqslant & \frac{\sigma^2\|c_1\|^2_2}{\lambda_{\min}(H(X)^T[\Phi_\rho(X,X)+\Sigma_{\gamma}]^{-1}H(X))}\nonumber\\
               \leqslant & \frac{\sigma^2\|c_1\|^2_2\lambda_{\max}(\Phi_\rho(X,X)+\Sigma_{\gamma})}{\lambda_{\min}(H(X)^TH(X))}.\label{eq:AppendixPartI}
\end{align}
Following development similar to above, $\lambda_{\min}(H(X)^TH(X))$ admits approximation
\begin{align}
\lambda_{\min}(H(X)^TH(X))=\lambda_{\min}(\sum^n_{i=1}h(x_i)h(x_i)^T)\approx m\lambda_{\min}(\int h(y)h(y)^TdF(y))=ms_2,\label{s2def}
\end{align}
with respect to the large sample distribution of the input locations, $F$.
Further, $s_2\geqslant 0$ with equality if and only if there exists $a\neq 0$ such that $h(y)^Ta=0$ with probability 1.
Combining (\ref{eq:AppendixPartII}) and (\ref{eq:AppendixPartI}) gives approximate upper bound for ${\rm Part(I)}+{\rm Part(II)}$
\begin{align}
& \frac{\sigma^2\|c_1\|^2_2\lambda_{\max}(\Phi_\rho(X,X)+\Sigma_{\gamma})}{ms_2} +\frac{\sigma^4\|c\|_2^2(m \sup_{u,v\in \Omega}\Phi_\rho(u,v)+\lambda_{{\rm max}}(\Sigma_{\gamma}))^2}{m^2s_1},\nonumber\\
                               \leqslant & \frac{\sigma^2\|c_1\|^2_2(m \sup_{u,v\in \Omega}\Phi_\rho(u,v)+\lambda_{{\rm max}}(\Sigma_{\gamma}))}{ms_2}+\frac{\sigma^4\|c\|_2^2(m \sup_{u,v\in \Omega}\Phi_\rho(u,v)+\lambda_{{\rm max}}(\Sigma_{\gamma}))^2}{m^2s_1},
\end{align}
finishing the proof of Theorem \ref{theopEstBound}.

\subsection{Proof of Proposition \ref{Uppc4}}\label{PfPropc4}
Recall that
\begin{align*}
(c_4)_t  =&\frac{\partial \hat{f}(x)}{\partial \tau_t}
     =\Phi_\rho(x,X)(\Phi_\rho(X,X)+\Sigma_{\gamma})^{-1}{\rm diag}(\frac{\partial \gamma_1}{\partial \tau_t}I_{k_1},...,\frac{\partial \gamma_i}{\partial \tau_t}I_{k_i},...,\frac{\partial \gamma_n}{\partial \tau_t}I_{k_n})\\
                &\quad\quad\quad\quad\times(\Phi_\rho(X,X)+\Sigma_{\gamma})^{-1}(y(X)-H(X)\beta).
\end{align*}
In this section we would give an upper bound of $(c_4)_t$. Without loss of generality, we can suppose $\Phi_\rho(x,x)=1$.
Let
\begin{align*}
\Phi_\rho(X,X)+\Sigma_{\gamma} =
\left[
\begin{array}{cc}
B_1+\Sigma_{\gamma_1} & R^T\\
R & B_2+\Sigma_{\gamma_2}
\end{array}\right],
\end{align*}
where
\begin{align*}
B_1 & = \mathbf{1}{\mathbf{1}^T},\\
\Sigma_{\gamma_1} & = \sigma_1^2I_{k_1},\\
R & = \Phi_\rho(X_2,x_1)\mathbf{1}^T,\\
B_2 & = \Phi_\rho(X_2,X_2).
\end{align*}
Thus, we have
\begin{align*}
 (\Phi_\rho(X,X)+\Sigma_{\gamma})^{-1}
 & =  \left[
\begin{array}{cc}
B_1+\Sigma_{\gamma_1} & R^T\\
R & B_2+\Sigma_{\gamma_2}
\end{array}\right]^{-1}\\
& = \left[
\begin{array}{cc}
B_{22}^{-1} & -B_{22}^{-1}R^T(B_2+\Sigma_{\gamma_2})^{-1}\\
-(B_2+\Sigma_{\gamma_2})^{-1}RB_{22}^{-1} & (B_2+\Sigma_{\gamma_2})^{-1}+(B_2+\Sigma_{\gamma_2})^{-1}RB_{22}^{-1}R^T(B_2+\Sigma_{\gamma_2})^{-1}
\end{array}\right],
\end{align*}
where $B_{22}=B_1+\Sigma_{\gamma_1}-R^T(B_2+\Sigma_{\gamma_2})^{-1}R$. Notice that
\begin{align*}
(B_1+\Sigma_{\gamma_1})^{-1}\mathbf{1}=\frac{1}{k_1+\sigma_1^2}\mathbf{1},
\end{align*}
we have
\begin{align*}
B_{22}^{-1} & = (B_1+\Sigma_{\gamma_1}-R^T(B_2+\Sigma_{\gamma_2})^{-1}R)^{-1}\\
            & = (B_1+\Sigma_{\gamma_1})^{-1}((B_1+\Sigma_{\gamma_1})^{-1}-(B_1+\Sigma_{\gamma_1})^{-1}R^T(B_2+\Sigma_{\gamma_2})^{-1}R(B_1+\Sigma_{\gamma_1})^{-1})^{-1}(B_1+\Sigma_{\gamma_1})^{-1}\\
            & = (B_1+\Sigma_{\gamma_1})^{-1}\bigg((B_1+\Sigma_{\gamma_1})^{-1}-\frac{1}{(k_1+\sigma_1^2)^2}\mathbf{1}\Phi_\rho(x_1,X_2)(B_2+\Sigma_{\gamma_2})^{-1}\Phi_\rho(X_2,x_1)\mathbf{1}^T\bigg)^{-1}(B_1+\Sigma_{\gamma_1})^{-1}\\
            & = (B_1+\Sigma_{\gamma_1})^{-1}\bigg((B_1+\Sigma_{\gamma_1})^{-1}-\frac{\Phi_\rho(x_1,X_2)(B_2+\Sigma_{\gamma_2})^{-1}\Phi_\rho(X_2,x_1)}{(k_1+\sigma_1^2)^2}\mathbf{1}\mathbf{1}^T\bigg)^{-1}(B_1+\Sigma_{\gamma_1})^{-1}.
\end{align*}
By binomial inverse theorem,
\begin{align*}
 & \bigg((B_1+\Sigma_{\gamma_1})^{-1}-\frac{\Phi_\rho(x_1,X_2)(B_2+\Sigma_{\gamma_2})^{-1}\Phi_\rho(X_2,x_1)}{(k_1+\sigma_1^2)^2}\mathbf{1}\mathbf{1}^T\bigg)^{-1}\\
& = B_1+\Sigma_{\gamma_1} + \frac{\Phi_\rho(x_1,X_2)(B_2+\Sigma_{\gamma_2})^{-1}\Phi_\rho(X_2,x_1)}{(k_1+\sigma_1^2)^2}\frac{(B_1+\Sigma_{\gamma_1})\mathbf{1}\mathbf{1}^T(B_1+\Sigma_{\gamma_1})}{1-\frac{\Phi_\rho(x_1,X_2)(B_2+\Sigma_{\gamma_2})^{-1}\Phi_\rho(X_2,x_1)}{(k_1+\sigma_1^2)^2}\mathbf{1}^T(B_1+\Sigma_{\gamma_1})\mathbf{1}}\\
& = B_1+\Sigma_{\gamma_1} + \frac{\Phi_\rho(x_1,X_2)(B_2+\Sigma_{\gamma_2})^{-1}\Phi_\rho(X_2,x_1)}{(k_1+\sigma_1^2)^2}\frac{(k_1+\sigma_1^2)^2\mathbf{1}\mathbf{1}^T}{1-\frac{\Phi_\rho(x_1,X_2)(B_2+\Sigma_{\gamma_2})^{-1}\Phi_\rho(X_2,x_1)}{(k_1+\sigma_1^2)^2}(k_1+\sigma_1^2)k_1}\\
& = B_1+\Sigma_{\gamma_1} + \Phi_\rho(x_1,X_2)(B_2+\Sigma_{\gamma_2})^{-1}\Phi_\rho(X_2,x_1)\frac{\mathbf{1}\mathbf{1}^T}{1-\frac{\Phi_\rho(x_1,X_2)(B_2+\Sigma_{\gamma_2})^{-1}\Phi_\rho(X_2,x_1)}{k_1+\sigma_1^2}k_1}.
\end{align*}
Let $d = \Phi_\rho(x_1,X_2)(B_2+\Sigma_{\gamma_2})^{-1}\Phi_\rho(X_2,x_1)$. Thus,
\begin{align*}
B_{22}^{-1} & =
(B_1+\Sigma_{\gamma_1})^{-1}\bigg(B_1+\Sigma_{\gamma_1} + d\frac{\mathbf{1}\mathbf{1}^T}{1-\frac{d}{k_1+\sigma_1^2}k_1}\bigg)(B_1+\Sigma_{\gamma_1})^{-1}\\
& = (B_1+\Sigma_{\gamma_1})^{-1} + d\frac{(B_1+\Sigma_{\gamma_1})^{-1}\mathbf{1}\mathbf{1}^T(B_1+\Sigma_{\gamma_1})^{-1}}{1-\frac{d}{k_1+\sigma_1^2}k_1}\\
& = (B_1+\Sigma_{\gamma_1})^{-1} + d\frac{\frac{1}{(k_1+\sigma_1^2)^2}\mathbf{1}\mathbf{1}^T}{1-\frac{d}{k_1+\sigma_1^2}k_1}.
\end{align*}
Since
\begin{align*}
\Phi_\rho(x,X) = (\Phi_\rho(x,x_1)\mathbf{1}^T,\Phi_\rho(x,X_2)),
\end{align*}
we have
\begin{align*}
&\Phi_\rho(x,X)(\Phi_\rho(X,X)+\Sigma_{\gamma})^{-1}\nonumber\\
 = &(\Phi_\rho(x,x_1)\mathbf{1}^T,\Phi_\rho(x,X_2))\\
& \times\left[
\begin{array}{cc}
B_{22}^{-1} & -B_{22}^{-1}R^T(B_2+\Sigma_{\gamma_2})^{-1}\\
-(B_2+\Sigma_{\gamma_2})^{-1}RB_{22}^{-1} & (B_2+\Sigma_{\gamma_2})^{-1}+(B_2+\Sigma_{\gamma_2})^{-1}RB_{22}^{-1}R^T(B_2+\Sigma_{\gamma_2})^{-1}
\end{array}\right],\\
 = &\bigg(
\Phi_\rho(x,x_1)\mathbf{1}^TB_{22}^{-1}-\Phi_\rho(x,X_2)(B_2+\Sigma_{\gamma_2})^{-1}RB_{22}^{-1},\\
& -\Phi_\rho(x,x_1)\mathbf{1}^TB_{22}^{-1}R^T(B_2+\Sigma_{\gamma_2})^{-1}+\Phi_\rho(x,X_2)(B_2+\Sigma_{\gamma_2})^{-1}\\
& +(B_2+\Sigma_{\gamma_2})^{-1}RB_{22}^{-1}R^T(B_2+\Sigma_{\gamma_2})^{-1}
\bigg).
\end{align*}
Notice that
\begin{align*}
B_{22}^{-1}\mathbf{1} & = (B_1+\Sigma_{\gamma_1})^{-1}\mathbf{1} + d\frac{\frac{1}{(k_1+\sigma_1^2)^2}\mathbf{1}\mathbf{1}^T\mathbf{1}}{1-\frac{d}{k_1+\sigma_1^2}k_1}\\
& = (B_1+\Sigma_{\gamma_1})^{-1}\mathbf{1} + d\frac{\frac{1}{(k_1+\sigma_1^2)^2}\mathbf{1}\mathbf{1}^T\mathbf{1}}{1-\frac{d}{k_1+\sigma_1^2}k_1}\\
& = (\frac{1}{k_1+\sigma_1^2} + \frac{\frac{dk_1}{(k_1+\sigma_1^2)^2}}{1-\frac{dk_1}{k_1+\sigma_1^2}})\mathbf{1},
\end{align*}
we have (let $d_1 = (\frac{1}{k_1+\sigma_1^2} + \frac{\frac{dk_1}{(k_1+\sigma_1^2)^2}}{1-\frac{dk_1}{k_1+\sigma_1^2}}) = \frac{1}{k_1+\sigma_1^2-dk_1}$)
\begin{align*}
& \Phi_\rho(x,x_1)\mathbf{1}^TB_{22}^{-1}-\Phi_\rho(x,X_2)(B_2+\Sigma_{\gamma_2})^{-1}RB_{22}^{-1}\\
= & \Phi_\rho(x,x_1)d_1\mathbf{1}^T-d_1\Phi_\rho(x,X_2)(B_2+\Sigma_{\gamma_2})^{-1}\Phi_\rho(X_2,x_1)\mathbf{1}^T,
\end{align*}
and
\begin{align*}
& -\Phi_\rho(x,x_1)\mathbf{1}^TB_{22}^{-1}R^T(B_2+\Sigma_{\gamma_2})^{-1}+\Phi_\rho(x,X_2)(B_2+\Sigma_{\gamma_2})^{-1}
 +(B_2+\Sigma_{\gamma_2})^{-1}RB_{22}^{-1}R^T(B_2+\Sigma_{\gamma_2})^{-1}\\
 = &  -\Phi_\rho(x,x_1)k_1d_1\Psi_\theta(x_1,X_2)(B_2+\Sigma_{\gamma_2})^{-1}+\Phi_\rho(x,X_2)(B_2+\Sigma_{\gamma_2})^{-1} \\
&  \quad + \Phi_\rho(x_1,X_2)(B_2+\Sigma_{\gamma_2})^{-1}\Phi_\rho(X_2,x_1)k_1d_1\Phi_\rho(x_1,X_2)(B_2+\Sigma_{\gamma_2})^{-1}\\
 = & -\Phi_\rho(x,x_1)k_1d_1\Phi_\rho(x_1,X_2)(B_2+\Sigma_{\gamma_2})^{-1}+\Phi_\rho(x,X_2)(B_2+\Sigma_{\gamma_2})^{-1} \\
& \quad +  dk_1d_1\Phi_\rho(x_1,X_2)(B_2+\Sigma_{\gamma_2})^{-1}.
\end{align*}
With the same procedure, we have
\begin{align*}
 & (\Phi_\rho(X,X)+\Sigma_{\gamma})^{-1}(y(X)-H(X)\beta)\\
 = & \bigg((y(x_1)-H(x_1)\beta)d_1{1}-d_1\Phi_\rho(x_1,X_2)(B_2+\Sigma_{\gamma_2})^{-1}(y(X_2)-H(X_2)\beta)\mathbf{1},\\
 & -(y(x_1)-H(x_1)\beta)k_1d_1(B_2+\Sigma_{\gamma_2})^{-1}\Phi_\rho(X_2,x_1)+(B_2+\Sigma_{\gamma_2})^{-1}(y(X_2)-H(X_2)\beta) \\
 + & dk_1d_1(B_2+\Sigma_{\gamma_2})^{-1}\Phi_\rho(X_2,x_1)\bigg).
\end{align*}
Thus,
\begin{align*}
(c_4)_t  =&\frac{\partial \hat{f}(x)}{\partial \tau_t}
     =\Phi_\rho(x,X)(\Phi_\rho(X,X)+\Sigma_{\gamma})^{-1}{\rm diag}(\frac{\partial \gamma_1}{\partial \tau_t}I_{k_1},...,\frac{\partial \gamma_i}{\partial \tau_t}I_{k_i},...,\frac{\partial \gamma_n}{\partial \tau_t}I_{k_n})\\
                &\quad\quad\quad\quad\times(\Phi_\rho(X,X)+\Sigma_{\gamma})^{-1}(y(X)-H(X)\beta)\\
= & k_1\frac{\partial \gamma_1}{\partial \tau_t}(\Phi_\rho(x,x_1)d_1-d_1\Phi_\rho(x,X_2)(B_2+\Sigma_{\gamma_2})^{-1}\Phi_\rho(X_2,x_1))((y(x_1)-H(x_1)\beta)d_1\\
& -d_1\Phi_\rho(x_1,X_2)(B_2+\Sigma_{\gamma_2})^{-1}(y(X_2)-H(X_2)\beta))\\
& +(-\Phi_\rho(x,x_1)k_1d_1\Phi_\rho(x_1,X_2)(B_2+\Sigma_{\gamma_2})^{-1}+\Phi_\rho(x,X_2)(B_2+\Sigma_{\gamma_2})^{-1} \\
 &  +  dk_1d_1\Phi_\rho(x_1,X_2)(B_2+\Sigma_{\gamma_2})^{-1}){\rm diag}(\frac{\partial \gamma_2}{\partial \tau_t}I_{k_2},...,\frac{\partial \gamma_i}{\partial \tau_t}I_{k_i},...,\frac{\partial \gamma_n}{\partial \tau_t}I_{k_n})\\
 & \times (-(y(x_1)-H(x_1)\beta)k_1d_1(B_2+\Sigma_{\gamma_2})^{-1}\Phi_\rho(X_2,x_1)+(B_2+\Sigma_{\gamma_2})^{-1}(y(X_2)-H(X_2)\beta) \\
 + & dk_1d_1(B_2+\Sigma_{\gamma_2})^{-1}\Phi_\rho(X_2,x_1))\\
\end{align*}
Let $d_2=\frac{1}{1+\sigma_1^2/k_1-d}$, we have
\begin{align*}
(c_4)_t  = & \frac{1}{k_1}\frac{\partial \gamma_1}{\partial \tau_t}(\Phi_\rho(x,x_1)d_2-d_2\Phi_\rho(x,X_2)(B_2+\Sigma_{\gamma_2})^{-1}\Phi_\rho(X_2,x_1))((y(x_1)-H(x_1)\beta)d_2\\
& -d_2\Phi_\rho(x_1,X_2)(B_2+\Sigma_{\gamma_2})^{-1}(y(X_2)-H(X_2)\beta))\\
& +(-\Phi_\rho(x,x_1)d_2\Phi_\rho(x_1,X_2)(B_2+\Sigma_{\gamma_2})^{-1}+\Phi_\rho(x,X_2)(B_2+\Sigma_{\gamma_2})^{-1} \\
 &  +  dd_2\Phi_\rho(x_1,X_2)(B_2+\Sigma_{\gamma_2})^{-1}){\rm diag}(\frac{\partial \gamma_2}{\partial \tau_t}I_{k_2},...,\frac{\partial \gamma_i}{\partial \tau_t}I_{k_i},...,\frac{\partial \gamma_n}{\partial \tau_t}I_{k_n})\\
 & \times (-(y(x_1)-H(x_1)\beta)d_2(B_2+\Sigma_{\gamma_2})^{-1}\Phi_\rho(X_2,x_1)+(B_2+\Sigma_{\gamma_2})^{-1}(y(X_2)-H(X_2)\beta) \\
 + & dd_2(B_2+\Sigma_{\epsilon_2})^{-1}\Phi_\rho(X_2,x_1))\\
 =&\Phi_\rho(x,X)(\Phi_\rho(X',X')+\Sigma_{\gamma})^{-1}{\rm diag}(\frac{1}{k_1}\frac{\partial \gamma_1}{\partial \tau_t},...,\frac{\partial \gamma_i}{\partial \tau_t}I_{k_i},...,\frac{\partial \gamma_n}{\partial \tau_t}I_{k_n})\\
                &\quad\quad\quad\quad\times(\Phi_\rho(X',X')+\Sigma_{\gamma})^{-1}(y(X')-H(X')\beta),
\end{align*}
where $X'=(x_1,X_2)$. Thus, by continuing this procedure, we have
\begin{align*}
|(c_4)_t| & \leqslant \frac{\|\Phi_\rho(x,\bar{X})\|_2\|\bar Y-H(\bar{X})\beta\|_2}{(\lambda_{{\rm min}}(\Phi_\rho(\bar{X},\bar{X})+\bar{\Sigma}_{\gamma}))^2}\max_{i:x_i\in\bar{X}}\bigg|\frac{1}{k_i}\frac{\partial \gamma_i}{\partial \tau_t}\bigg|.
\end{align*}

\subsection{Proof of Theorem \ref{Thm:ParaEstNum}}\label{pfParaEstnum}

The following lemma, which describes the accuracy of solving linear systems \citep{golub}, will be used to develop a bound on the numeric error.
\begin{lemma}\label{lemmaLinearSystems}
Suppose ${Ax}={b}$ and $\tilde{{A}}\tilde{{x}}=\tilde{{b}}$ with $\|\tilde{{A}}-{A}\|_2\leqslant \delta \|{A}\|_2$, $\|\tilde{{b}}-{b}\|_2\leqslant \delta \|{b}\|_2$, and $\kappa({A})=r/\delta<1/\delta$ for some $\delta>0$. Then, $\tilde{{A}}$ is non-singular,
\begin{align}
&\frac{\|\tilde{{x}}\|_2}{\|{x}\|_2}\leqslant \frac{1+r}{1-r},\nonumber\\
&\frac{\|\tilde{{x}}-{x}\|_2}{\|{x}\|_2}\leqslant \frac{2\delta}{1-r}\kappa({A}),\label{eq:NumericLinearSystem}
\end{align}
where $\kappa({A})=\|{A}\|_2\|{A}^{-1}\|_2$.
\end{lemma}
Further, for conformable ${A}$, ${b}$, $\tilde{{A}}$, and $\tilde{{b}}$, we have
\begin{align}
& \|{A}{b}-\tilde{{A}}\tilde{{b}}\|_2=\|{A}({b}-\tilde{{b}})-(\tilde{{A}}-{A})
\tilde{{b}}\|_2\nonumber\\
\leqslant & \|{A}({b}-\tilde{{b}})\|_2+\|(\tilde{{A}}-{A})\tilde{{b}}\|_2
\leqslant \|{A}\|_2\|({b}-\tilde{{b}})\|_2+\|(\tilde{{A}}-{A})\|_2\|\tilde{{b}}\|_2.\label{eq:NumericMinus}
\end{align}

In order to satisfy the conditions of Lemma \ref{lemmaLinearSystems}, we make a few assumptions \emph{in addition to Assumption \ref{numericAssump}}, in particular, with regard to the accuracy of numeric optimization.
\begin{assumption}\label{paramNumericAssump}
Assume $\kappa(\hat{A})=r/\delta$ with $r<1$ and
\begin{align*}
\|\hat{A}-\tilde{A}\|_2\leqslant \delta\|\hat{A}\|_2,\|\Psi_{\hat{\theta}}(\bar X,x)-\Psi_{\tilde{\theta}}(\bar X,x)\|_2\leqslant \delta\|\Psi_{\hat{\theta}}(\bar X,x)\|_2.
\end{align*}
\end{assumption}
Note that this assumption does not concern the parameter estimates themselves, but instead the accuracy of the solution to the optimization problem.
If the optimization problem is solved with sufficient accuracy, then this assumption will be satisfied.
However, as we will see in the following, the regression function coefficients $\beta$ have great potential to cause problems. Briefly, in order to control parameter estimation numeric error, we need that numeric properties are even more tightly controlled,
in particular, an even smaller condition number of $\Psi_{\theta}(\bar X,\bar X)+\Sigma_{\epsilon}$, which is stated in the following assumption.
\begin{assumption}\label{assumpInParaNum}
\begin{align*}
\delta\kappa(\hat{A})\kappa(H(\bar X)^TH(\bar X))\bigg(1+(1+\delta)^2+\frac{(1+\delta)^2}{1-r}\kappa(\hat{A})\bigg)<1.
\end{align*}
\end{assumption}
Assumption \ref{assumpInParaNum} is a strong assumption, since it requires $\delta\kappa(\hat{A})^2$ to be relatively small, at least smaller than 1. However, since our goal is to make $\kappa(\hat{A})$ small, in practice this condition is not too difficult to achieve, since we can control the condition number of $\hat{A}$.\par
The following lemma states that if Assumption \ref{assumpInParaNum} holds, combining Assumption \ref{paramNumericAssump}, the conditions of Lemma \ref{lemmaLinearSystems} holds.
\begin{lemma}\label{lemmacondofPnEst}
Let
\begin{align}
r_1&=\delta\kappa(\hat{A})\kappa(H(\bar X)^TH(\bar X))\bigg(1+(1+\delta)^2+\frac{(1+\delta)^2}{1-r}\kappa(\hat{A})\bigg)\nonumber\\
&\quad+\frac{1}{2}\min\{\delta,1-\delta\kappa(\hat{A})\kappa(H(\bar X)^TH(\bar X))\bigg(1+(1+\delta)^2+\frac{(1+\delta)^2}{1-r}\kappa(\hat{A})\bigg)\}\nonumber\\
\delta_1 & = \frac{r_1}{\kappa(H(\bar X)^T\hat{A}^{-1}H(\bar X))}.\label{conditionofrdelta}
\end{align}
Suppose Assumptions \ref{numericAssump}, \ref{paramNumericAssump}, and \ref{assumpInParaNum} hold, we have $r_1<1$ and
\begin{align}
 \|H(\bar X)^T\hat{A}^{-1}H(\bar X)-\tilde{H}(\bar X)^T\tilde{A}^{-1}\tilde{H}(\bar X)\|_2 <\delta_1\|H(\bar X)^T\hat{A}^{-1}H(\bar X)\|_2. \label{eq:ParaNumericHX}
\end{align}
\end{lemma}
Thus, we have all tools to give an upper bound of $|\hat{f}_{\hat{\vartheta}}-\hat{f}_{\tilde{\vartheta}}|$. By triangle inequality,
\begin{align}
&|\hat{f}_{\hat{\vartheta}}-\hat{f}_{\tilde{\vartheta}}|\nonumber\\
& =|h(x)^T\hat{\beta}+\Psi_{\hat{\theta}}(x,\bar X)\hat{A}^{-1}(\bar Y-H(\bar X)\hat{\beta})
                                                        -(h(x)^T\tilde{\beta}+\Psi_{\tilde{\theta}}(x,\bar X)\tilde{A}^{-1}(\bar Y - H(\bar X)\tilde{\beta}))|\nonumber\\
& = |h(x)^T(\hat{\beta}-\tilde{\beta})+\bar Y^T(\hat{A}^{-1}\Psi_{\hat{\theta}}(\bar X,x)-\tilde{A}^{-1}\Psi_{\tilde{\theta}}(\bar X,x))\nonumber\\
& -[\Psi_{\hat{\theta}}(x,\bar X)\hat{A}^{-1}H(\bar X)\hat{\beta}-\Psi_{\tilde{\theta}}(x,\bar X)\tilde{A}^{-1}H(\bar X)\tilde{\beta}]|\nonumber\\
& \leqslant \|h(x)\|_2\|\hat{\beta}-\tilde{\beta}\|_2+\|\bar Y\|_2\|\hat{A}^{-1}\Psi_{\hat{\theta}}(\bar X,x)-\tilde{A}^{-1}\Psi_{\tilde{\theta}}(\bar X,x)\|_2\nonumber\\
& + \|\Psi_{\hat{\theta}}(x,\bar X)\hat{A}^{-1}H(\bar X)\hat{\beta}-\Psi_{\tilde{\theta}}(x,\bar X)\tilde{A}^{-1}H(\bar X)\tilde{\beta}\|_2\nonumber\\
& = {\rm Part}(i)+{\rm Part}(ii)+{\rm Part}(iii).\label{eq:ParaNumericTotal1}
\end{align}
Part$(ii)$ can be bounded using Lemma \ref{lemmaLinearSystems} as
\begin{align}
{\rm Part}(ii) & \leqslant \|\bar Y\|_2\frac{2\delta}{1-r}\kappa(\hat{A})\|\hat{A}^{-1}\Psi_{\hat{\theta}}(\bar X,x)\|_2.\label{eq:ParaNumericTotal2}
\end{align}
Similarly, Part$(iii)$ can be bounded using (\ref{eq:NumericMinus}) and Lemma \ref{lemmaLinearSystems} as
\begin{align}
&{\rm Part}(iii) \leqslant \|H(\bar X)\|_2\|\hat{\beta}\|_2\frac{2\delta}{1-r}\kappa(\hat{A})\|\hat{A}^{-1}\Psi_{\hat{\theta}}(\bar X,x)\|_2+
\|H(\bar X)\|_2\|\hat{\beta}-\tilde{\beta}\|_2\|\tilde{A}^{-1}\Psi_{\tilde{\theta}}(\bar X,x)\|_2\nonumber\\
& \leqslant \|H(\bar X)\|_2\|\hat{\beta}\|_2\frac{2\delta}{1-r}\kappa(\hat{A})\|\hat{A}^{-1}\Psi_{\hat{\theta}}(\bar X,x)\|_2+
\|H(\bar X)\|_2\|\hat{\beta}-\tilde{\beta}\|_2\frac{1+r}{1-r}\|\hat{A}^{-1}\Psi_{\hat{\theta}}(\bar X,x)\|_2.\label{eq:ParaNumericTotal3}
\end{align}
Combining (\ref{eq:ParaNumericTotal1}), (\ref{eq:ParaNumericTotal2}) and (\ref{eq:ParaNumericTotal3})
gives
\begin{align}
|\hat{f}_{\hat{\vartheta}}-\hat{f}_{\tilde{\vartheta}}|&\leqslant \frac{2\delta\kappa(\hat{A})}{(1-r)\lambda_{\min}(\hat{A})} \|\Psi_{\hat{\theta}}(\bar X,x)\|_2(\|\bar Y\|_2+\|H(\bar X)\|_2\|\hat{\beta}\|_2)\nonumber\\
&\quad +\|\hat{\beta}-\tilde{\beta}\|_2(\|h(x)\|_2+\frac{1+r}{(1-r)\lambda_{\min}(\hat{A})}\|H(\bar X)\|_2\|\Psi_{\hat{\theta}}(\bar X,x)\|_2).\label{eq:ParaNumericTotal}
\end{align}
Notice that the first term in (\ref{eq:ParaNumericTotal}) can be controlled by restraining $g(\Sigma_M,\Sigma_{\epsilon})$, as defined in (\ref{eq:Numericg}). The second part can be controlled by, in addition, restraining $\|\hat{\beta}-\tilde{\beta}\|_2$.
Recall that
\begin{align*}
\hat{\beta}=(H(\bar X)^T\hat{A}^{-1}H(\bar X))^{-1}H(\bar X)^T\hat{A}^{-1}\bar Y,\\
\tilde{\beta}=(\tilde{H}(\bar X)^T\tilde{A}^{-1}\tilde{H}(\bar X))^{-1}\tilde{H}(\bar X)^T\tilde{A}^{-1}\tilde{Y}.
\end{align*}

Since by Lemma \ref{lemmacondofPnEst}, the condition of Lemma \ref{lemmaLinearSystems} holds. Thus, by Lemma \ref{lemmaLinearSystems}, we have
\begin{align*}
\|\hat{\beta}-\tilde{\beta}\|_2 & \leqslant \frac{2\delta_1}{1-r_1}\kappa(H(\bar X)^T\hat{A}^{-1}H(\bar X))\|\|\hat{\beta}\|_2 =\frac{2r_1}{1-r_1}   \|\hat{\beta}\|_2.
\end{align*}
By plugging in (\ref{conditionofrdelta}), we have
\begin{align}
\|\hat{\beta}-\tilde{\beta}\|_2 & \leqslant 2\delta\bigg(\kappa(\hat{A})\kappa(H(\bar X)^TH(\bar X))\bigg(1+(1+\delta)^2+\frac{(1+\delta)^2}{1-r}\kappa(\hat{A})\bigg)+1\bigg)\|\hat{\beta}\|_2.\label{eq:ParaBetaUpper}
\end{align}

Combining (\ref{eq:ParaNumericTotal}) and ($\ref{eq:ParaBetaUpper}$), we finish the proof.

\subsection{Proof of Lemma \ref{lemmacondofPnEst}}
Notice that if Assumption \ref{assumpInParaNum} holds, we have $r_1<1$. We only need to prove
(\ref{eq:ParaNumericHX}). Notice that
\begin{align}
& \|H(\bar X)^T\hat{A}^{-1}H(\bar X)-\tilde{H}(\bar X)^T\tilde{A}^{-1}\tilde{H}(\bar X)\|_2\nonumber\\
\leqslant & \|H(\bar X)^T\hat{A}^{-1}H(\bar X)-\tilde{H}(\bar X)^T\hat{A}^{-1}H(\bar X)\|_2+\|\tilde{H}(\bar X)^T\hat{A}^{-1}H(\bar X)-\tilde{H}(\bar X)^T\tilde{A}^{-1}\tilde{H}(\bar X)\|_2\nonumber\\
\leqslant & \delta \|H(\bar X)\|_2\|\hat{A}^{-1}H(\bar X)\|_2+\|\tilde{H}(\bar X)^T\hat{A}^{-1}H(\bar X)-\tilde{H}(\bar X)^T\tilde{A}^{-1}\tilde{H}(\bar X)\|_2\nonumber\\
\leqslant & \delta \|H(\bar X)\|^2_2\|\hat{A}^{-1}\|_2+\|\tilde{H}(\bar X)^T\hat{A}^{-1}H(\bar X)-\tilde{H}(\bar X)^T\tilde{A}^{-1}\tilde{H}(\bar X)\|_2,\label{eq:ParaNumericHX1}
\end{align}
where the first inequality is true because of the triangle inequality, the second inequality is true because of Assumption \ref{paramNumericAssump}, and the third inequality is true because $\|G^{-1}d\|_2\leqslant \|G^{-1}\|_2\|d\|$ for any vector $d$ and non-singular matrix $G$.
The second term in (\ref{eq:ParaNumericHX1}) has
\begin{align}
& \|\tilde{H}(\bar X)^T\hat{A}^{-1}H(\bar X)-\tilde{H}(\bar X)^T\tilde{A}^{-1}\tilde{H}(\bar X)\|_2\nonumber\\
\leqslant & \|\tilde{H}(\bar X)^T\hat{A}^{-1}H(\bar X)-\tilde{H}(\bar X)^T\hat{A}^{-1}\tilde{H}(\bar X)\|_2+\|\tilde{H}(\bar X)^T\hat{A}^{-1}\tilde{H}(\bar X)-\tilde{H}(\bar X)^T\tilde{A}^{-1}\tilde{H}(\bar X)\|_2\nonumber\\
\leqslant & \delta \|\tilde{H}(\bar X)\|_2\|\hat{A}^{-1}\tilde{H}(\bar X)\|_2+\|\tilde{H}(\bar X)^T(\hat{A}^{-1}-\tilde{A}^{-1})\tilde{H}(\bar X)\|_2\nonumber\\
\leqslant & \delta \|\tilde{H}(\bar X)\|^2_2\|\hat{A}^{-1}\|_2+\|\hat{A}^{-1}-\tilde{A}^{-1}\|_2\|\tilde{H}(\bar X)\|^2_2\nonumber\\
\leqslant & \delta (1+\delta)^2\|H(\bar X)\|^2_2\|\hat{A}^{-1}\|_2+(1+\delta)^2\|H(\bar X)\|^2_2\|\hat{A}^{-1}-\tilde{A}^{-1}\|_2,\label{eq:ParaNumericHX2}
\end{align}
where the first inequality is true is because of the triangle inequality, the second inequality is true because of Assumption \ref{paramNumericAssump}, the third inequality is true because $\|G^{-1}d\|_2\leqslant \|G^{-1}\|_2\|d\|$, and the last inequality is true because by Assumption \ref{numericAssump}, $\|\tilde{H}(\bar X)\|_2\leqslant (1+\delta)\|H(\bar X)\|_2$.
Next, $\|\hat{A}^{-1}-\tilde{A}^{-1}\|_2$ is bounded.


For any $x\in \mathbb{R}^n$ such that $\|x\|_2=1$, let $y_1,y_2\in \mathbb{R}^n$ such that $\hat{A}y_1=x$ and $\tilde{A}y_2=x$. Let $\delta_A=\tilde{A}-\hat{A}$. Thus, $(\hat{A}+\delta_A) y_2=x$. Notice that by assumption,
\begin{align*}
\|\hat{A}^{-1}\delta_A\|_2\leqslant \delta\|\hat{A}^{-1}\|_2\|\hat{A}\|_2=r<1\quad{\rm and}\quad
(I+\hat{A}^{-1}\delta_A)y_2=y_1.
\end{align*}
The following Lemma from \cite{golub} will be used.
\begin{lemma}\label{inverseLemma}
Suppose $F\in \mathbb{R}^{n\times n}$, $\|F\|_2<1$. Then $I-F$ is invertible and
\begin{align*}
\|(I-F)^{-1}\|_2\leqslant \frac{1}{1-\|F\|_2},
\end{align*}
where $I$ is identity matrix in $\mathbb{R}^{n\times n}$.
\end{lemma}
By Lemma \ref{inverseLemma}, we have
\begin{align*}
\|y_2\|_2\leqslant \|(I+\hat{A}^{-1}\delta_A)^{-1}\|_2\|y_2\|_2\leqslant \frac{1}{1-r}\|y_1\|_2\quad{\rm and}\quad
y_1-y_2=\hat{A}^{-1}\delta_Ay_2.
\end{align*}
So,
\begin{align*}
\|y_1-y_2\|_2 & \leqslant \|\hat{A}^{-1}\delta_A\|_2\|y_2\|_2\leqslant \delta\|\hat{A}^{-1}\|_2\|\hat{A}\|_2\|y_2\|_2=\frac{\delta}{1-r}\kappa(\hat{A})\|y_1\|_2.
\end{align*}
Plugging in $y_1$ and $y_2$ gives
\begin{align}
\|(\hat{A}^{-1}-\tilde{A}^{-1})x\|_2 & \leqslant \frac{\delta}{1-r}\kappa(\hat{A})\|\hat{A}^{-1}x\|_2\leqslant \frac{\delta}{1-r}\kappa(\hat{A})\|\hat{A}^{-1}\|_2= \frac{\delta}{1-r}\kappa(\hat{A})\frac{1}{\lambda_{{\rm min}}(\hat{A})},\label{eq:ParaNumericLin}
\end{align}
indicating
\begin{align}
\|\hat{A}^{-1}-\tilde{A}^{-1}\|_2 & \leqslant\frac{\delta}{1-r}\kappa(\hat{A})\frac{1}{\lambda_{{\rm min}}(\hat{A})},\label{eq:ParaNumericA}
\end{align}
since (\ref{eq:ParaNumericLin}) is true for any $x$ with $\|x\|_2=1$.
Combining (\ref{eq:ParaNumericHX1}), (\ref{eq:ParaNumericHX2}), and (\ref{eq:ParaNumericA}) gives
\begin{align}
& \|H(\bar X)^T\hat{A}^{-1}H(\bar X)-\tilde{H}(\bar X)^T\tilde{A}^{-1}\tilde{H}(\bar X)\|_2\nonumber\\
\leqslant & \delta \|H(\bar X)\|^2_2\|\hat{A}^{-1}\|_2+\delta (1+\delta)^2\|H(\bar X)\|^2_2\|\hat{A}^{-1}\|_2+(1+\delta)^2\|H(\bar X)\|^2_2\|\hat{A}^{-1}-\tilde{A}^{-1}\|_2\nonumber\\
\leqslant & \delta \frac{\|H(\bar X)\|^2_2}{\lambda_{{\rm min}}(\hat{A})}+\delta (1+\delta)^2\frac{\|H(\bar X)\|^2_2}{\lambda_{{\rm min}}(\hat{A})}+\frac{\delta(1+\delta)^2}{1-r}\kappa(\hat{A})\frac{\|H(\bar X)\|^2_2}{\lambda_{{\rm min}}(\hat{A})}\nonumber\\
= & \frac{\delta\|H(\bar X)\|^2_2}{\lambda_{{\rm min}}(\hat{A})}\bigg(1+(1+\delta)^2+\frac{(1+\delta)^2}{1-r}\kappa(\hat{A})\bigg).\label{eq:ParaNumericHXTotal}
\end{align}
Thus, (\ref{eq:ParaNumericHX}) holds if
\begin{align*}
\frac{\delta\|H(\bar X)\|^2_2}{\lambda_{{\rm min}}(\hat{A})}\bigg(1+(1+\delta)^2+\frac{(1+\delta)^2}{1-r}\kappa(\hat{A})\bigg)<\delta_1\|H(\bar X)^T\hat{A}^{-1}H(\bar X)\|_2,
\end{align*}
or equivalently
\begin{align}
&\frac{\delta\|H(\bar X)\|^2_2}{\lambda_{{\rm min}}(\hat{A})\lambda_{{\rm min}}(H(\bar X)^T\hat{A}^{-1}H(\bar X))}\bigg(1+(1+\delta)^2+\frac{(1+\delta)^2}{1-r}\kappa(\hat{A})\bigg)<\delta_1\kappa(H(\bar X)^T\hat{A}^{-1}H(\bar X)).\label{eq:ParaNumericCond1}
\end{align}
Next, we simplify (\ref{eq:ParaNumericCond1}). Notice that the left-hand side of (\ref{eq:ParaNumericCond1}) has
\begin{align*}
& \frac{\delta\|H(\bar X)\|^2_2}{\lambda_{{\rm min}}(\hat{A})\lambda_{{\rm min}}(H(\bar X)^T\hat{A}^{-1}H\bar (X))}\bigg(1+(1+\delta)^2+\frac{(1+\delta)^2}{1-r}\kappa(\hat{A})\bigg)\nonumber\\
\leqslant & \frac{\lambda_{{\rm max}}(\hat{A})}{\lambda_{{\rm min}}(H(\bar X)^TH(\bar X))}\frac{\delta\|H(\bar X)\|^2_2}{\lambda_{{\rm min}}(\hat{A})}\bigg(1+(1+\delta)^2+\frac{(1+\delta)^2}{1-r}\kappa(\hat{A})\bigg)\nonumber\\
= & \delta\kappa(\hat{A})\kappa(H(\bar X)^TH(\bar X))\bigg(1+(1+\delta)^2+\frac{(1+\delta)^2}{1-r}\kappa(\hat{A})\bigg),
\end{align*}
so, if
\begin{align}
\delta\kappa(\hat{A})\kappa(H(\bar X)^TH(\bar X))\bigg(1+(1+\delta)^2+\frac{(1+\delta)^2}{1-r}\kappa(\hat{A})\bigg)<\delta_1\kappa(H(\bar X)^T\hat{A}^{-1}H(\bar X)),\label{eq:conditionOfHX}
\end{align}
(\ref{eq:ParaNumericHX}) holds. By plugging in (\ref{conditionofrdelta}), we have (\ref{eq:conditionOfHX}) holds, which finishes the proof.

\bibliography{papref}

\end{document}